%
%

\magnification=\magstep1
\documentstyle{amsppt}

\parskip=\baselineskip 
\baselineskip=1.1\baselineskip
\parindent=0pt
\loadmsbm
\UseAMSsymbols
\raggedbottom

\def\pt{\hbox{\it pt}}
\def\Vol{\hbox{vol}}
\def\sol{\operatorname{sol}}
\def\vor{\operatorname{vor}}
\def\dih{\operatorname{dih}}
\def\arccot{\operatorname{arccot}}
\def\score{\sigma}
\def\doct{\delta_{oct}}
\def\dihmin{\dih_{\min}}
\def\dihmax{\dih_{\max}}
\def\rad{\operatorname{rad}}

\def\qed{{\hbox{}\nobreak\hfill\vrule height8pt width6pt depth 0pt}\medskip}
\def\diag|#1|#2|{\vbox to #1in {\vskip.3in\centerline{\tt Diagram #2}\vss} }
\def\v{\hskip -3.5pt }
\def\gram|#1|#2|#3|{
        {
        \smallskip
        \hbox to \hsize
        {\hfill
        \vrule \vbox{ \hrule \vskip 6pt \centerline{\it Diagram #2}
         \vskip #1in %
             \includegraphics{#3}\hrule }
        \v\vrule\hfill
        }
\smallskip}}

\centerline{\bf Sphere Packings I}
\bigskip
\centerline{Thomas C. Hales}
\bigskip

\footnote""{
I would like to thank D. J. Muder
for the appendix and the referees for suggesting other
substantial improvements.\hskip -13pt
}
\footnote""{\line{\it\hfill published in Discrete and Computational Geometry, 17:1-51, 1997}}

{\bf Abstract:}  We describe a program to prove the Kepler
conjecture on sphere packings.  We then carry out the first
step of this program.  Each packing determines a 
decomposition of space into Delaunay simplices, which
are grouped together into finite configurations
called Delaunay stars.  A score, which is related to the
density of packings,  is assigned to each Delaunay
star.  We conjecture that the score of every
Delaunay star is at most the score of
the stars in the face-centered cubic and hexagonal
close packings.  This conjecture implies the Kepler conjecture.
To complete the first step of the program, we show
that every Delaunay star that satisfies a certain regularity
condition satisfies the conjecture.

\bigskip

{\bf Contents:}
1. Introduction, 2. The Program,
3. Quasi-regular Tetrahedra, 4. Quadrilaterals, 5. Restrictions,
6. Combinatorics, 
7. The Method of Subdivision, 8. Explicit Formulas for
Compression, Volume, and Angle, 9. Floating-Point Calculations.

Appendix. D. J. Muder's Proof of Theorem 6.1.

\bigskip
\centerline{\bf Section 1. Introduction}
\bigskip

The Kepler conjecture
asserts that no packing of
spheres in three dimensions has density exceeding that of the
face-centered cubic lattice packing. This density is 
$\pi/\sqrt{18}\approx 0.74048$.
  In an earlier paper \cite{H2}, we showed how to
reduce the Kep\-ler conjecture to a finite calculation.
That paper also gave numerical evidence in support
of the method and conjecture.
This finite calculation is a series of optimization problems
involving up to 53 spheres in an explicit compact region of
Euclidean space.
Computers have little difficulty in treating problems of this size
numerically, but a naive attempt to make a thorough study
of the possible arrangements of these spheres would quickly
exhaust the world's computer resources.

The first purpose of this paper is to describe a
program designed to give a rigorous proof of the Kepler conjecture.
A sketch of a related program appears in \cite{H2}.  Although
the approach of \cite{H2} is based on substantial
numerical evidence,
some of the constructions of that paper are needlessly complicated.
This paper streamlines some of those constructions and
replaces others with constructions that are more amenable
to rigorous methods.
For this program to succeed, the
original optimization problem must be partitioned into a series
of much smaller problems that may be treated by current computer
technology or hand calculation.

The second purpose of this paper is to carry out the first step of
the proposed program.  A statement of the result is
contained in Theorem 1 below.  

Background to another approach to this problem is found in 
\cite{H3}.
To add more detail to the proposed program, 
we recall some constructions from earlier papers \cite{H1}, \cite{H2}.
Begin with
a packing of nonoverlapping spheres of radius 1 in Euclidean three-space.
The {\it density} of a packing is defined in \cite{H1}.  It is defined
as a limit of the ratio of the
volume of the unit balls in a large region of space
to the volume of the large region.
The density of the packing may be improved by adding spheres until
there is no further room to do so.  The resulting packing is said to
be {\it saturated\/}.  It has the property that no point in space has
distance greater than $2$ from the center of some sphere.

Every saturated packing gives rise to a decomposition of space into
simplices called the {\it Delaunay decomposition}.  The vertices
of each Delaunay simplex are centers of spheres of the packing.  
None of the centers
of the spheres of the packing lie in the interior of the circumscribing
sphere of any Delaunay simplex.  
In fact, this property is enough to
completely determine the Delaunay decomposition except for
certain degenerate packings.  A degeneracy occurs,
for instance, when two Delaunay simplices have the same circumscribing
sphere.  In practice, these degeneracies are important, because they
occur in the face-centered cubic and hexagonal close packings.
The paper \cite{H2} shows how to resolve the degeneracies by taking a
small perturbation of the packing.  In general, the
Delaunay decomposition will depend on this perturbation.
We refer to the centers of the packing as {\it vertices\/},
since the structure of the simplicial
decomposition of space will be our primary concern.
For a proof that the Delaunay decomposition is a dissection of space
into simplices, we refer the reader to \cite{R}.

The Delaunay decomposition is dual to the well-known Voronoi
decomposition.  
If the vertices of the Delaunay simplices are in nondegenerate
position, two vertices are joined by an edge exactly when
the two corresponding Voronoi cells share a face, three vertices
form a face exactly when the three Voronoi cells share an edge, and
four vertices form a simplex exactly when the four corresponding
Voronoi cells share a vertex.  In other words, two vertices are
joined by an edge if they lie on a sphere that does not contain
any other of the vertices, and so forth (again assuming the vertices
to be in nondegenerate position).  The collection of all simplices that
share a given vertex is called a {\it Delaunay star}.  (This is
a provisional definition: it will
be refined below.)

Every Delaunay simplex has
edges between 2 and 4 in length and, because of the saturation
of the packing,  a circumradius of at
most 2. 
We assume that every simplex $S$ in this paper comes with a fixed
order on its edges, $1,\ldots,6$. 
The order on the edges is to be arranged so that
the first, second, and third edges meet at a vertex.  We may
also assume that the edges numbered $i$ and $i+3$ are opposite
edges for $i=1,2,3$.  We define $S(y_1,\ldots,y_6)$ to be
the (ordered) simplex whose $i$th edge has length $y_i$.
If $S$ is a Delaunay simplex in a fixed Delaunay star, then
it has a distinguished vertex, the vertex common to all simplices
in the star.  In this situation, we assume that the edges
are numbered so that the first, second, and third edges meet
at the distinguished vertex.

A function, known as the {\it compression\/} $\Gamma(S)$, is defined on
the space of all Delaunay simplices.  Let 
$\delta_{oct} = (-3\pi+12\arccos(1/\sqrt{3}))/\sqrt{8}\approx 0.720903$
be the density of a regular octahedron with edges of length 2.
That is, place a unit ball at each vertex of the
octahedron, and let $\doct$ be the ratio of the volume of the
part of the balls in the octahedron to the volume of the octahedron.
Let $S$ be a Delaunay simplex.  Let $B$ be the union of
four unit balls placed at each of the vertices of $S$.
Define the compression as
$$\Gamma(S) = -\delta_{oct} \Vol(S) + \Vol(S\cap B).$$
We extend the definition of compression to Delaunay stars $D^*$
by setting $\Gamma(D^*) = \sum\Gamma(S)$, with the sum running
over all the Delaunay simplices in the star.  

In this and subsequent work, we single out for
special treatment the edges of length between 2 and 2.51.
The constant $2.51$ was determined experimentally to have a number of
desirable properties.  This constant will appear throughout
the paper.  We will call vertices that come within $2.51$ of each
other {\it close neighbors}.

We say that the convex hull of four vertices
is a 
{\it quasi-regular tetrahedron\/} (or simply a {\it tetrahedron})
if all four vertices are close neighbors of one another.
Suppose that we have a configuration of
six vertices in bijection with the vertices
of an octahedron with the property that two vertices are close
neighbors if and only if the corresponding vertices of the
octahedron are adjacent.
Suppose further that exactly one of the three diagonals has
length at most $2\sqrt{2}$.
  In this
case we call the convex hull of the six vertices a
{\it quasi-regular
octahedron\/} (or simply an {\it octahedron}).

The compatibility of quasi-regular tetrahedra and octahedra with 
the Delaunay decomposition
is established in Section 3.
We think of Euclidean space as
the union of quasi-regular tetrahedra, octahedra, and various
less interesting Delaunay simplices.  From now on, a Delaunay
star is to be the collection of all quasi-regular tetrahedra, octahedra,
and Delaunay simplices 
that share a common vertex $v$.  This collection of
Delaunay simplices and quasi-regular solids is often, but not
always, the same as the objects called Delaunay stars in \cite{H2}.  We
warn the reader of this shift in terminology.

It is convenient to measure the compression in multiples of
the compression of the regular simplex of edge length 2.
We define a {\it point} (abbreviated $\pt$) to be
$\Gamma(S(2,2,2,2,2,2))$. 
We have $\pt = 
        {11\pi/3 -12\arccos(1/\sqrt{3})}\approx 0.0553736$.

One of the main purposes
of this paper and its sequel is to replace the compression by
a function (called the {\it score}) that has better properties
than the compression.  
Further details on the definition of score will appear in Section 2.
We are now able to state the main theorem of this paper. 

\bigskip
{\bf Theorem 1.}  If a Delaunay star is composed entirely of quasi-regular
tetrahedra, then its score is less than $8\,\pt$.
\bigskip

The idea of the proof is the following.  Consider the unit sphere
whose center is the center 
of the Delaunay star $D^*$.  The intersection of a simplex
in $D^*$ with this unit sphere is a spherical triangle.  
For example,
a regular tetrahedron with edges of length 2
gives a triangle on the unit sphere of arc
length $\pi/3$.  
The star
$D^*$ gives a triangulation of the unit sphere.  
The restriction
on the lengths of the edges of a quasi-regular tetrahedron
constrains the triangles
in the triangulation.
We classify all triangulations that potentially come from a star
scoring more than $8\,\pt$.
Section 5 develops
a long list of properties that must be satisfied by 
the triangulation of a high-scoring star.

It is then necessary to classify all the triangulations that possess
the properties on this list.
The original classification was carried
out by a computer program, which generated all
potential triangulations and checked them against the list.
D. J. Muder has made a significant improvement in the argument by giving
a direct, computer-free classification.  
His result appears in the appendix.

As it turns out, there is only one triangulation that satisfies all
of the properties on the list.  Section 7 proves that 
Delaunay stars with this triangulation
score less than $8\,\pt$.  This will complete the main
thread of the argument.

There are a number of
estimates in this paper that are established by computer.  These
estimates are used throughout the paper, even though their proofs
are not discussed until Sections 8 and 9.
These sections may be viewed as a series of technical appendices
giving explicit formulas for the compression, dihedral angles,
solid angles, volumes, and other quantities that must be estimated.
The final section states the inequalities and gives
details about the computerized verification.
There is no vicious
circle here:  the results of Sections 8 and 9 do not
rely on anything from Sections 2--7.

\bigskip

There are several functions of a Delaunay simplex that will be
used throughout this paper.  The compression $\Gamma(S)$ has
been defined above.  The {\it dihedral angle\/} $\dih(S)$ is defined
to be the dihedral angle of the simplex $S$ along the first
edge (with respect to the fixed order on the edges of $S$).
Set 
$\dihmin =  \dih(S(2,2,2.51,2,2,2.51))\approx 0.8639$;
$\dihmax = \dih(S(2.51,2,2,2.51,2,2)) = 
	\arccos(-29003/96999) \approx 1.874444$.
We will see that $\dihmin$ and $\dihmax$ are lower and upper bounds
on the dihedral angles of quasi-regular tetrahedra.
The {\it solid angle\/} (measured in steradians) at the vertex
joining the first, second, and third edges is denoted
$\sol(S)$.  
The intersection of $S$ with the ball of unit
radius centered at this vertex has volume $\sol(S)/3$
(see \cite{H1,2.1}).  For example, $\sol(S(2,2,2,2,2,2))\approx 0.55$.
Let $\rad(S)$ be the circumradius of the simplex $S$.
In Section 2, we will define two other functions:
$\vor(S)$, which is related to the volume of Voronoi cells,
and the score $\score(S)$.
Finally, let $\eta(a,b,c)$ denote the circumradius of a triangle with
edges $a$, $b$, $c$.
Explicit formulas for all these functions appear in Section 8.

\bigskip
\centerline{\bf Section 2. The  Program}

By proving Theorem 1, the main purpose of this paper will be
achieved.  Nevertheless, it might be helpful to give a series of
comments about how Theorem 1 may be viewed as the
solution to the first of a handful of optimization problems that
would collectively provide a solution to the Kepler conjecture.

We begin with some notation and terminology.  We fix
a Delaunay star $D^*$ about a vertex $v_0$, which we
take to be the origin, and we consider the unit sphere at $v_0$.
Let $v_1$ and  $v_2$ be vertices of $D^*$ such that
$v_0$, $v_1$, and $v_2$ are all close neighbors of one another.
We take the radial projections $p_i$
of $v_i$ to the unit sphere with center at the origin
and connect the points $p_1$ and $p_2$ by a geodesic
arc on the sphere.  
We mark all such arcs on the unit sphere.
Lemma 3.10 will show that the arcs meet only at their endpoints.
The closures of the
connected components of the complement of these arcs are
regions on the unit sphere,
called the {\it standard regions}.  We may remove the arcs
that do not bound one of the regions.
The resulting system of edges and regions will be
referred to as the {\it standard decomposition\/} of the unit sphere.

Let $C$ be the cone with vertex $v_0$ over one of the 
standard regions.
The collection of the Delaunay simplices,
quasi-regular tetrahedra, and quasi-regular octahedra
of $D^*$ in $C$ (together with the distinguished vertex $v_0$)
will be called a {\it standard cluster}.  Each
Delaunay simplex in $D^*$ belongs to a unique standard cluster.
Each triangle in the standard decomposition of the unit sphere
is associated
with a unique quasi-regular tetrahedron, and each tetrahedron
determines a triangle in the standard decomposition
(Lemma 3.7). We
may identify quasi-regular tetrahedra with clusters over
triangular regions.

We assign a score to each standard cluster in \cite{H4,3}.  
In this
section we define the score of a quasi-regular tetrahedron and
describe the properties that the score should have in general.

Let $S$ be a quasi-regular tetrahedron.  It
is a standard cluster in a Delaunay star with center $v_0$.
If the circumradius
of $S$ is at most $1.41$, then we define the score to be
$\Gamma(S)$.  

If the circumradius is greater than $1.41$,
then embed the simplex $S$ in Euclidean three-space.
Partition Euclidean space into four infinite regions (infinite
Voronoi cells) by associating with each vertex of $S$ the points of
space closest to that vertex.  By intersecting $S$ with each
of the four regions, we partition $S$ into four pieces $\hat S_0$,
$\hat S_1$, $\hat S_2$, and $\hat S_3$, 
corresponding to its four vertices $v_0$,
$v_1$, $v_2$, and $v_3$. Let $\sol_i$ be the solid angle at the
vertex $v_i$ of the simplex.
The expression $-4\doct\Vol(\hat S_0) + 4\sol_0/3$ 
is an analytic function of the lengths of the edges for
simplices $S$ that contain their circumcenters.  
(Explicit formulas appear in Section 8.)  This
function may be analytically continued to a function
of the lengths of the edges for simplices $S$ that do
not necessarily contain their circumcenters.
Let $\vor(S)$ be defined as the analytic continuation
of  $-4\doct \Vol(\hat S_0) + 4\sol_0/3$.  

In this case, 
define the score of $S$ to be $\vor(S)$.  Write $\score(S)$
for the score of a quasi-regular tetrahedron.
In summary, the score is
$$\score(S) = \cases \Gamma(S), & 
	\text{if the circumradius of $S$ is at most $1.41$,}\\
	\vor(S), & \text{otherwise.}\endcases
$$

The quasi-regular tetrahedron $S$ appears in four Delaunay stars.  In
the other three Delaunay stars, the distinguished vertices will
be $v_1$, $v_2$, and $v_3$, so that $S$, viewed as a standard cluster
in the other Delaunay stars, will have scores 
$-\doct\Vol(\hat S_i) + \sol_i/3$
(or their analytic continuations),
for $i=1,2,3$.  
By definition,
$$\Gamma(S) = 
	\sum_{i=0}^3 (-\doct \Vol(\hat S_i)+\sol_i/3).$$
The sum of the scores of $S$, for each
of the four vertices of $S$, is
$4\Gamma(S)$.  This is the same total that is obtained
by summing the compression of $S$ at each of its vertices.
This is the property we need to relate the score to
the density of the packing.  It means that although the
score reapportions the compression among neighboring Delaunay
stars, the average of the compression over a large region
of space equals the average of the score over the same region, 
up to a negligible boundary term.

The analytic continuation 
in the definition of $\vor(S)$
has the following geometric interpretation.
If the circumcenter of a quasi-regular tetrahedron
$S$ is not contained in $S$, a small tip of the infinite
Voronoi cell at $v_0$ (or some other vertex) will protrude through
the opposite face of the Delaunay simplex. The volume of this
 small protruding tip
is not counted in $-4\doct\Vol(\hat S_0) + 4\sol_0/3$, but it is
counted in the analytic continuation.  The analytic continuations
of the scores of $S$ for each of the other three vertices
acquire a term representing the negative volume of a part
of the tip. The three parts together constitute the entire
tip, so that the negative volumes exactly offset the
volume of the tip, and the sum of the four scores
of $S$ is $4\Gamma(S)$.  Details appear in \cite{H4}.

The general definition of the score will have similar properties.
To each standard cluster of a Delaunay star $D^*$ a score will be 
assigned.
The rough idea is to
let the score of a simplex in a cluster be 
the compression $\Gamma(S)$
if the circumradius of every face of $S$ is small, and otherwise to
let the score be defined by Voronoi cells (in a way
that generalizes the definition for quasi-regular tetrahedra).

The score $\sigma(D^*)$
of a Delaunay star is defined as the sum of
the scores of its standard clusters.  The score has
the following properties \cite{H4,3.1 and 3.5}. 

1.  The score of a standard cluster depends only on 
	the cluster, and not on the way it sits in a Delaunay
	star or in the Delaunay decomposition of space.

2.  The Delaunay stars of the face-centered cubic and hexagonal
	close packings score exactly $8\,\pt$.

3.  The score is asymptotic to the compression over large
regions of space.  We make this more precise.  Let $\Lambda$
denote the vertices of a saturated packing.  Let $\Lambda_N$
denote the vertices inside the ball of radius $N$. (Fix
a center.)  Let $D^*(v)$ denote the Delaunay star at
$v\in \Lambda$.  Then the score satisfies (in Landau's
notation)
$$\sum_{\Lambda_N} \score(D^*(v)) = \sum_{\Lambda_N} \Gamma(D^*(v)) + O(N
^2).$$

{\bf Lemma 2.1.}  If the score of every Delaunay star
in a saturated packing is at most $s<16\pi/3$, 
then the density
of the packing is at most $16\pi\doct/(16\pi-3s)$.
If the score of every Delaunay star in a
packing is at most $8\,\pt$, then the density of the
packing is at most $\pi/\sqrt{18}$.

{\bf Proof:}  The second claim is the special case
$s=8\,\pt$.  The proof relies on property 3.
The number of vertices such that $D^*(v)$ meets the boundary
of the ball $B_N$ of radius $N$ has order $O(N^2)$.
Since the Delaunay stars give a four-fold cover of ${\Bbb R}^3$,
we have
$$\align
4(-\doct\Vol(B_N)+|\Lambda_N|{4\pi\over3}) &=
  \sum_{\Lambda_N}(-\doct\Vol(D^*(v))+4\sol(D^*(v))/3) + O(N^2)\\
  &=\sum_{\Lambda_N}\Gamma(D^*(v))+O(N^2)\\
  &=\sum_{\Lambda_N}\sigma(D^*(v)) +O(N^2)\\
  &\le s|\Lambda_N| + O(N^2).
\endalign
$$
Rearranging, we get
$${4\pi|\Lambda_N|\over3\Vol(B_N)}\le
	{\doct\over(1- 3s/16\pi)} + {O(N^2)\over \Vol(B_N)}.$$
In the limit, the left-hand side is the density and the right-hand
side is the bound.  Similar arguments 
can be found in \cite{H1} and \cite{H2}.\qed

\bigskip

The following conjecture is fundamental.  By the lemma, 
this conjecture implies the Kepler conjecture.  The lemma
also shows that weaker bounds than $8\,\pt$
on the score might be used to give new upper bounds
on the density of sphere packings.

\bigskip
{\bf Conjecture 2.2.}  The score of every Delaunay star is at most $8\,\pt$.

The basic philosophy behind the approach of this paper is that quasi-regular
tetrahedra are the only clusters that give a positive score,  standard clusters
over quadrilateral regions should be the
only other clusters that 
may give a score of zero, and every other standard cluster
should give a negative score.  Moreover, we will prove that no
quasi-regular tetrahedron gives more than $1\,\pt$.

Thus, heuristically, we try to obtain a high score by including
as many triangular regions as possible.  If we allow any other
shape, preference should be given to  quadrilaterals.  Any other
shape of region should be avoided if possible.  If these other
regions occur, they should be accompanied by additional
triangular regions to compensate for the negative score of the
region.  We will see later in this paper that even triangular
regions tend to give a low score unless they are arranged to give five
triangles around each vertex.

The main steps in a proof of the Kepler conjecture are

{\it

\def\ha{ \hangindent=20pt \hangafter=1\relax }
1. A proof that even if all 
regions are triangular the total score
is less than $8\,\pt$

\ha
2. A proof that standard clusters in regions of more than
three sides score at most $0\,\pt$

\ha
3. A proof that if all of the 
standard regions 
are triangles or quadrilaterals,
then the total score is less than $8\,\pt$ (excluding the
case of pentagonal prisms)

\ha
4.  A proof that if some standard region 
has more than
four sides, then the star scores less than $8\,\pt$

\ha
5.  A proof that pentagonal prisms score less than $8\,\pt$

}

The division of the problem into these steps is
quite arbitrary.  
  They were originally
intended to be steps of roughly equal magnitude, although
is has turned out that a construction in
\cite{H4} has made the second step substantially easier
than the long calculations of the third step.

This paper carries out the first step.   The second step of the
program is also complete \cite{H4}.
Partial results are known for the third step \cite{H5}.
In the fourth step,
it will be necessary to argue
that these regions take up too much space, give too little in return,
and have such strongly incompatible shapes that they cannot be
part of a winning strategy.

To make step 5 precise,
we define {\it pentagonal prisms\/} to be Delaunay stars whose
standard decomposition has ten triangles and five quadrilaterals,
with the five quadrilaterals in a band around the equator,
capped on both ends by five triangles
(Diagram 2.3.)
The conjecture in this section
asserts, in particular,
that pentagonal prisms, which created such difficulties
in \cite{H2},  score less than $8\,\pt$.  
The final step has been separated from the third
step, because the estimates are expected to be more delicate
for pentagonal prisms than for a general
Delaunay star in the third step.

\gram|2.2|2.3|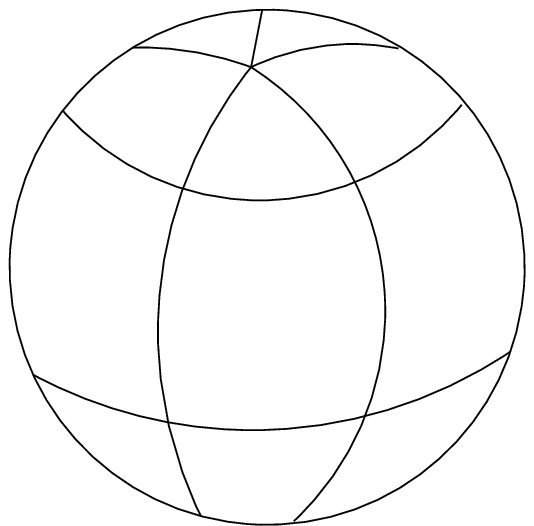|  

One of the main shortcomings of the compression is that 
pentagonal prisms have compression
greater than $8\,\pt$ (see \cite{H2}).
Numerical evidence suggests that
the upper bound on the compression is attained
by a pentagonal prism, denoted $D^*_{\operatorname{ppdp}}$ in \cite{H2},
at about $8.156\,\pt$,
and this means that the link between the compression
and the Kepler conjecture is indirect.
The score appears to correct this shortcoming.  

What evidence is there for the conjecture and program?
They come as the result of extensive computer experimentation.
I have checked the conjecture against much of the data
obtained in the numerical studies of \cite{H2,9.3}.  
The data suggest that the score tends to give a dramatic
improvement over the compression, often improving
the bound by more than a point.  The score
of the particularly troublesome 
 pentagonal prism $D^*_{\operatorname{ppdp}}$ 
drops safely under $8\,\pt$.  
I have checked a broad assortment
of other pentagonal prisms and have found them all to score less
than $8\,\pt$.  

The second step shows that no serious pathologies can arise.
The only way to form a Delaunay star with a positive score
is by arranging a number of quasi-regular tetrahedra
around a vertex (together with other standard clusters
that can only lower the score).  There must be at least
eight tetrahedra to score $8\,\pt$, and if there are any
distortions in these tetrahedra, there must be at least
nine.  However, as this paper shows, too many quasi-regular
tetrahedra in any star are also harmful.  Future papers
will impose additional limits on the structure of the
optimal Delaunay star.

\bigskip
\centerline{\bf Section 3.  Quasi-Regular Tetrahedra}
\bigskip

This section studies the compatibility of the Delaunay simplices
and the quasi-regular solids.
Fix three
vertices $v_1$, $v_2$, and $v_3$ that are close neighbors to
one another.
Let $T$
be the triangle with vertices $v_i$.
It does not follow
that $T$ is the face of a Delaunay simplex.
However,
as we will see, when $T$ is not the face of a simplex, the
arrangement of the surrounding simplices is almost completely determined.

If $T$ is not the face of a Delaunay simplex, then
we will show that there are two additional vertices
$v_0$ and $v_0'$, where $v_0$ and $v_0'$ are close neighbors
to $v_1$, $v_2$, and $v_3$.  This means that there
are two quasi-regular tetrahedra $S_1$ and $S_2$ with
vertices $(v_0,v_1,v_2,v_3)$ and $(v_0',v_1,v_2,v_3)$, 
respectively,
that have the common face $T$ (see Diagram 3.1.a).
We will see that $S_1\cup S_2$ is the union of three
Delaunay simplices with vertices $(v_0,v_0',v_1,v_2)$, 
$(v_0,v_0',v_2,v_3)$,
and $(v_0,v_0',v_3,v_1)$ 
in Diagram 3.1.b.  This section establishes
that this is the only situation in which quasi-regular
tetrahedra are not Delaunay simplices:  they must come
in pairs and their union must be three Delaunay simplices
joined along a common edge.
The decomposition of this paper is obtained by
taking each such triple of Delaunay simplices (3.1.b)
and replacing the triple by a pair of quasi-regular
tetrahedra (3.1.a).

\gram|2.2|3.1|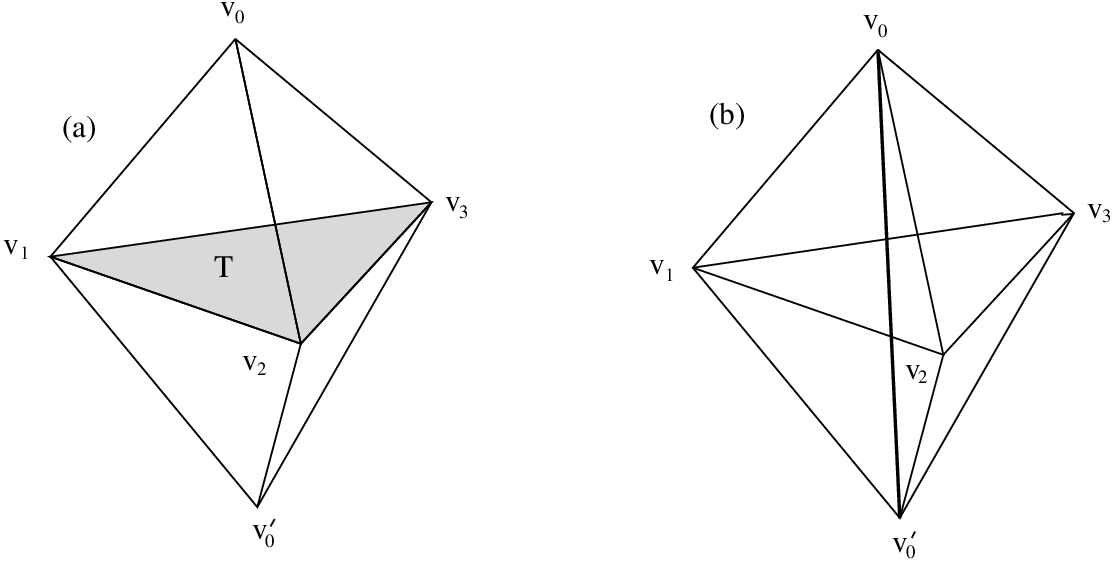|  

From the dual perspective of Voronoi cells, the Voronoi
cell at $v_0$ (or $v_0'$)  would have a small tip protruding from $S_1$
through $T$, if the vertex $v_0'$ were not present.
The vertex $v_0'$ slices off this protruding tip so that
the Voronoi cells at $v_0$ and $v_0'$ have a small face in common.

\bigskip 
{\bf Lemma 3.2.}  Suppose that the circumradius of the triangle $T$
is less than $\sqrt{2}$.  Then $T$ is the face of a Delaunay simplex.

\bigskip
{\bf Proof:}  Let $r<\sqrt{2}$ be the radius of the circle that
circumscribes $T$, and let $c$ be the center of the circle.
The sphere of radius $r$ at $c$ does not contain any vertices of $D^*$
other than $v_1$, $v_2$, and $v_3$.  By the definition of the
Delaunay decomposition (as described in
the introduction), this implies that $T$
is the face of a simplex. \qed

\bigskip
{\bf Remark 3.3.} We have several constraints on the edge lengths, if
$T$ is not a face of a Delaunay simplex.
Consider the circumradius $\eta(a,b,c)$ of 
a triangle whose edges have lengths $a,b,c$
between $2$ and $2.51$.
Since $2.51^2 < 2^2 + 2^2$, we see that the
triangle is acute, so that $\eta(a,b,c)$ 
is monotonically increasing in $a$, $b$, and $c$.
This
gives simple estimates relating the circumradius to $a$, $b$, and $c$.
The circumradius is at most 
$\eta(2.51,2.51,2.51)=2.51/\sqrt{3}\approx 1.449$.
If the circumradius is at least $\sqrt{2}\approx 1.41421$, then
$a$, $b$, and $c$ are greater than $2.3$  
($\eta(2.3,2.51,2.51)\approx 1.41191<\sqrt{2}$).
Under the same hypothesis, two of $a$, $b$, and $c$
are greater than $2.41$  
($\eta(2.41,2.41,2.51)<\sqrt{2}$).  Finally,
at least one edge has length greater than
$2.44$
($\eta(2.44,2.44,2.44)<\sqrt{2}$).

\bigskip
{\bf Lemma 3.4.}  Let $T$ be the triangle with vertices
$v_i$.  Assume that the vertices $v_i$ are close
neighbors of one another.
Suppose there is a vertex $v_0$ that lies 
closer to the circumcenter
of $T$ than the vertices of $T$ do.
Then the vertex $v_0$ satisfies $2\le |v_0-v_i|< 2.15$,
for $i=1,2,3$.  In particular, the convex hull of $v_0,\ldots,v_3$
is a quasi-regular
tetrahedron $S$ with face $T$.

\bigskip
Another way of stating the hypothesis on the circumcenter
is to say that the plane of $T$ separates $v_0$ from the
circumcenter of $S$.
Because of the constraints on the edge lengths in Remark 3.3,
the three other
faces of $S$ are faces of Delaunay simplices.

\bigskip
{\bf Proof:}  We defer the proof to Section 8.2.5.
\qed

\bigskip
{\bf Lemma 3.5.}  Let $v,v_1,v_2,v_3$, and $v_4$ be distinct
vertices with pairwise distances at least $2$.  Suppose that
the pairs $(v_i,v_j)$ are close neighbors for $\{i,j\}\ne\{1,4\}$.
Then $v$ does not lies in the convex hull of $(v_1,v_2,v_3,v_4)$.

\bigskip
{\bf Proof:} For a contradiction, suppose $v$ lies in the convex hull.
Since $|v-v_i| \ge 2$ and
for $\{i,j\}\ne\{1,4\}$,
 $v_i$ and $v_j$ are close neighbors,
the angle formed
by $v_i$ and $v_j$ at the vertex $v$ is at most
$\theta_0 = \arccos(1-2.51^2/8) \approx 1.357$.
Each such pair of vertices gives a geodesic arc of length at most
$\theta_0$ radians on the unit sphere centered at $v$.
We obtain in this
way a triangulation of the unit sphere by 
four triangles,
two with edges of length
at most $\theta_0=2\arcsin(2.51/4)$ radians, and two others
that fit together to form a quadrilateral with edges of at most $\theta_0$ radians.
By the spherical law of cosines, the area formula for a spherical triangle,
and \cite{H2,6.1}, each of the first two triangles has 
area at most
$3\arccos(\cos\theta_0/(1+\cos\theta_0)) -\pi \approx 1.04$.  
By the same lemma,
the quadrilateral
has area at most the area of a regular quadrilateral of side $\theta_0$, or
about $2.8$.
Since the combined area of the two triangles
and quadrilateral is less than $4\pi$, they cannot give the
desired triangulation.
To see that $v$ cannot lie on the boundary, it is enough to check
that a triangle having two edges of lengths between $2$ and $2.51$
cannot contain a point that has distance $2$ or more from each
vertex.  We leave this as an exercise.\qed

\bigskip

{\bf Corollary 3.5.1}  No vertex of the packing is ever an interior point 
of a quasi-regular tetrahedron or octahedron.

\bigskip
{\bf Proof:}  The corollary is immediate for a tetrahedron.
For an octahedron, draw the distinguished diagonal 
and apply the lemma to
each of the four resulting simplices.\qed
\bigskip

Let $T$ be a triangle of circumradius between
$\sqrt{2}$ and $2.51/\sqrt{3}$,
with edges of length between $2$ and $2.51$.
Consider the line through the circumcenter of $T$, perpendicular
to the plane of $T$.   Let $s$ be the finite segment in this
line whose endpoints are the circumcenters of the two simplices
with face $T$ formed by placing an additional vertex at distance
two from the three vertices of of $T$ on either side
of the plane through $T$.

\bigskip
{\bf Lemma 3.6.} Let $S$ be a simplex formed by the vertices
of $T$ and a fourth vertex $v_0'$.  Suppose that the circumcenter
of $S$ lies on the segment $s$.  Assume that $v_0'$ has distance
at least $2$ from each of 
the vertices of $T$.  Then $v'_0$ has distance
less than 2.3 from each of the vertices of $T$.

{\bf Proof:}  Let $v_1$, $v_2$, $v_3$ be the
vertices of $T$. For a contradiction, assume that 
$v_0'$ has distance at least $2.3$ from a vertex $v_1$ of $T$.
Lemma 3.4 shows that the plane through $T$ does not separate $v_0'$
from the circumcenter of $S$.   If the circumcenter lies in $s$
(and is not separated from $v_0'$
by the plane through $T$), then by moving $v_0'$ to
decrease the circumradius, the circumcenter remains in $s$.

Let $S$ be the simplex with vertices $v_0'$, $v_1$, $v_2$,
and $v_3$.
The circumcenter of $S$ lies in the interior of $S$.  We omit
the proof, because it is established by methods similar to
(but longer than)
the proof of Lemma 3.4.  
Thus,
the circumradius is increasing in the lengths of $|v_0'-v_i|$,
for $i=1,2,3$ (see 8.2.4).  

Let $R$ be the circumradius of a simplex with face $T$
and center an endpoint of $s$.
 We will prove that the circumradius of $S$
is greater than
$R$, contrary to our hypothesis.
Moving $v'_0$ to decrease 
the circumradius,
we may take the distances to $v_i$
to be precisely $2.3$, $2$, and $2$.
We may move $v_1$, $v_2$, and $v_3$ along their fixed circumscribing circle
until $|v_2-v_3|=2.51$ and $|v_1-v_2|=|v_1-v_3|$ in a way
that does not decrease any of the distances from $v_0'$ to $v_i$.
Repeating the previous step, we may retain our assumption that $v_0'$ has
distances exactly $2.3$, $2$, and $2$ from the vertices $v_i$
as before.  We have reduced
the problem to a one-dimensional family of tetrahedra parametrized
by the radius $r$ of the circumscribing circle of $T$.  Set
$x(r) = |v_1-v_3|=|v_2-v_3|$.  To obtain our desired contradiction,
we must show that the circumradius $R'(r)$ of the simplex
$S(2.3,2,2,2.51,x(r),x(r))$ satisfies $R'(r)>R(r)$.  Since
both $R'$ and $R$ are increasing in $r$, for $r\in[\sqrt{2},2.51/\sqrt{3}]$,
the desired inequality follows if we evaluate the 200 constants
$$R'(r_i) - R(r_{i+1}),\qquad \hbox{for } r_i = \sqrt{2} +
                \left({2.51\over\sqrt{3}}-\sqrt{2}\right){i\over 200},$$
for $i=0,\ldots,199$, and check that they are all positive.
(The smallest is about $0.00005799$, which occurs for $i=199$.)\qed

Let $T$ be a triangle made up of three close neighbors.  Suppose
that $T$ is not the face of a Delaunay simplex.  
There exists a vertex $v_0$ whose distance to the circumcenter
of $T$ is less than the circumradius
of $T$.  Let $S$ be the quasi-regular
tetrahedron
formed by $v_0$ and the vertices of $T$.  It is not a Delaunay
simplex, so there exists
a vertex $v_0'$ that is less than the circumradius of $S$ from the
circumcenter of $S$.  Let $S'$ be the simplex formed
by $v_0'$ and the vertices of $T$.  It is not a Delaunay simplex
either.  The circumcenter of $S$ lies in the
segment $s$ of Lemma 3.4, so the circumcenter of $S'$ does too.
The lengths of the edges of $S$ and $S'$ other than $T$ are
constrained by Lemmas 3.4 and 3.6.  In particular, in light
of Remark 3.3, the faces other than $T$ of $S$ and $S'$ are
faces of Delaunay simplices.

If $v_0$ and $v_0'$ lie on the same side of the plane through $T$,
then either $v_0'$ lies in $S$, $v_0$ lies in $S'$, or the faces
of $S$ and $S'$ intersect nonsimplicially.  None of these
situations can occur because a nondegenerate 
Delaunay decomposition is
a Euclidean simplicial complex and because of Lemma 3.5.
We conclude that $v_0$ and $v_0'$ lie on opposite sides of the plane 
through  $T$.

$S\cup S'$ is bounded by Delaunay faces, so $S\cup S'$ is a union
of Delaunay simplices.  The fourth vertex of the
Delaunay simplex in $S\cup S'$ with face $(v_0,v_1,v_2)$ cannot
be $v_3$ ($S$ is not a Delaunay simplex), so
it must be $v_0'$.
Similarly, $(v_0,v_0',v_2,v_3)$ and $(v_0,v_0',v_3,v_1)$ are
Delaunay simplices.  These three Delaunay simplices cannot be
quasi-regular tetrahedra by Lemma 8.3.2. 

The assumption that $T$ is not a face has completely determined
the surrounding geometry:  there are two quasi-regular tetrahedra
$S$ and $S'$ along $T$ such that $S\cup S'$ is a union of three
Delaunay simplices.

\bigskip
{\bf Lemma 3.7.}  
Let $L$ be a union of standard
regions.  Suppose that the boundary of
$L$ consists of three edges.  Then either $L$ or its complement is a single 
triangle.

\smallskip

For example, the interior of $L$ cannot have the form
of Diagram 3.8.  Lemma 3.5 (proof) shows that if all regions are
triangles, then there are at least $12$ triangles, so that
the exterior of $L$ cannot have the form of Diagram 3.8 either.

\gram|2|3.8|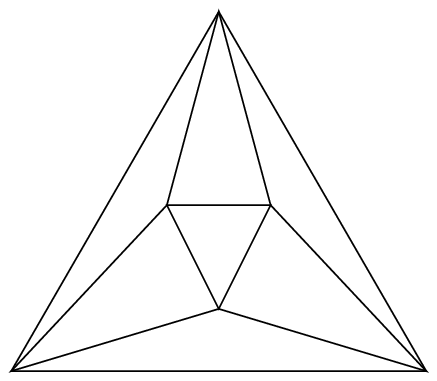|   

{\bf Proof:} Replacing $L$ by its complement if necessary, we may assume
that the area of $L$ is less than its complement.  The triangular boundary
corresponds to four vertices $v_0$ (the origin), $v_1$, $v_2$, and $v_3$.
The close-neighbor constraints on the lengths 
show that
the convex hull of $v_0,\ldots,v_3$ is a quasi-regular tetrahedron.
By construction each quasi-regular tetrahedron is a single cluster.
\qed

We say that a point $v\in{\Bbb R}^3$ is {\it enclosed\/} by a region
on the unit sphere if the interior of the cone 
(with vertex $v_0$) over that region contains
$v$.  For example, in Diagram 3.9, the point $v$ is enclosed by
the given spherical triangle.  

\gram|2.1|3.9|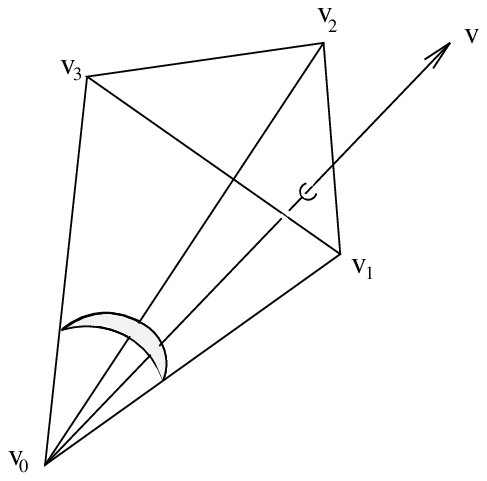|  

\smallskip

The following lemma was used in Section 2 to define the standard
decomposition.

{\bf Lemma 3.10.}  Fix a Delaunay star $D^*$ with center $v_0$.
Draw geodesic arcs on the unit sphere at $v_0$ for every
triple of close neighbors $v_0$, $v_1$, $v_2$ (as in Section 2).
The resulting system of arcs do not meet except at endpoints.

{\bf Proof.}  
Our proof will be based on the fact that a nondegenerate
Delaunay decomposition is a Euclidean simplicial complex.
 Let $T_1$ and $T_2$
be two triangles made from two such triples of close neighbors.
We have $T_i\subset S_i\cup S_i'$, where $S_i$ and $S_i$ are
the Delaunay simplices with face $T_i$ if $T_i$ is the face
of a Delaunay simplex, and they are the two quasi-regular
tetrahedra with face $T_i$ constructed above, otherwise.
Since a nondegenerate Delaunay decomposition is a 
Euclidean simplicial complex, $S_1\cup S_1'$
meets $S_2\cup S_2'$ simplicially.  By the restrictions on the
lengths of the edges in Lemmas 3.4 and 3.6, this forces 
$T_1$ to
 intersect $T_2$ simplicially.
The result follows.\qed

\bigskip\hbox{}\bigskip
\centerline{\bf Section 4. Quadrilaterals}

Fix a Delaunay star composed entirely of quasi-regular tetrahedra and
consider the associated triangulation of the unit sphere.  Let $L$ be
a region of the sphere bounded by four edges of the triangulation.
$L$ will be the union of two or more triangles.  Replacing $L$ by
its complement in the unit sphere if necessary, we assume that the area
of $L$ is less than that of its complement.

\bigskip
We claim that in this context, $L$ is the union of either
two or four triangles, as illustrated in Diagram 4.1.  In particular,
$L$ encloses at most one vertex.
If a diagonal to the quadrilateral $L$ is an edge of the
triangulation, the region $L$ is divided into two triangles each
associated with a quasi-regular tetrahedron.  In particular, there
is no enclosed vertex.  (Lemma 3.7
precludes any subtriangulation of a triangular region.)  If, however,
there is a single enclosed vertex and neither
diagonal is an edge of the triangulation,  then the only possible
triangulation of $L$ is the one of the diagram.  Proposition 4.2
completes the proof of the claim.

\bigskip
\gram|1.4|4.1|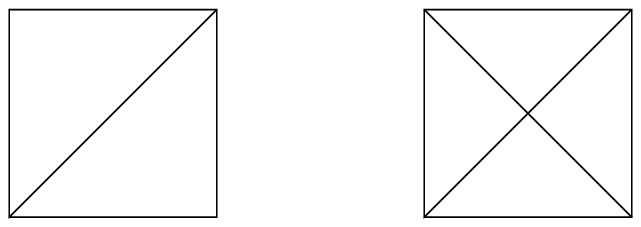|
\bigskip

{\bf Proposition 4.2.}  A union of regions
(of area less than $2\pi$)
 bounded by exactly four
edges 
cannot enclose two vertices of distance at most $2.51$ from the origin.

\bigskip
  This argument
is somewhat delicate:  if our parameter 2.51 had been set at 2.541,
for instance, such an arrangement would exist.  First we prove a useful
reduction.

\bigskip
{\bf Lemma 4.3.}  Assume a figure exists with vectors $v_1,\ldots,v_4$,
$v$, and $v'$ subject to the constraints
$$\align
2\le &|v_i|\le 2.51,\\
2\le&|v_i-v_{i+1}|\le k_i, \\
2\le&|v_i-v_{i+2}|,\\
2\le&|v-v'|,\\
h_i\le &|w-v_i|, \\
2\le &|w|\le \ell, \hbox{ for } w=v,v' \hbox{ and }
        i=1,\ldots,4 \ (\hbox{mod} 4)\\
\endalign
$$
where $\ell$, $h_i$, and $k_i$ are fixed constants that satisfy
$\ell\in[2.51,2\sqrt{2}]$,
$h_i\in[2,2\sqrt{2})$, $k_i\in[2,2.51]$.
Let $L$ be the quadrilateral on the unit sphere with vertices $v_i/|v_i|$
and edges running between consecutive vertices.  Assume that $v$ and $v'$
lie in the cone at the origin obtained by scaling $L$.
Then another figure exists made of a (new) collection of vectors
$v_1,\ldots,v_4$, $v$, and $v'$ subject to the constraints
above together with the additional constraints
$$\align
&|v_i-v_{i+1}|=k_i\\
&|v_i|=2, \hbox { for } i=1,\ldots,4,\\
&|v|=|v'|=\ell.
\endalign
$$
Moreover, the quadrilateral $L$ may be assumed to be convex.
\bigskip

{\bf Proof (4.3):}
By rescaling $v$ and $v'$, we may assume that $|v|=|v'|=\ell$.
(Moving $v$ or $v'$ away
from the origin increases its distance from the other vertices of the
configuration.)

The diagonals satisfy $|v_2-v_4|>2.1$ and $|v_1-v_3|>2.1$.  
Otherwise,
if say $|v_1-v_3|\le2.1$, then the faces with vertices $(v_1,v_2,v_3)$
and $(v_1,v_4,v_3)$ have circumradius less than $\sqrt{2}$.
The edge from $0$ to $v$ has length at most $\sqrt{2}$, so this edge
cannot intersect these faces by the Euclidean simplicial complex
argument used before 3.7 and in 3.10.  By Lemma 3.5, $v$
cannot lie in the convex hull of $(0,v_1,v_3,v_i)$.  This leaves
$v$ nowhere to go, and a figure with $|v_1-v_3|\le2.1$ does not exist.

Next, we claim that we may assume that the quadrilateral $L$
is convex (in the sense that it contains the geodesic arcs
between any two points in the region).  To see this, suppose the vertex
$v_i$
lies in the cone over the convex hull of the other three vertices $v_j$.
  Consider
the plane $P$ through the origin, $v_{i-1}$, and $v_{i+1}$.  The reflection
$v'_i$ of $v_i$ through $P$ is no closer to $v$, $v'$, or
$v_{i+2}$  and has the same distance
to the origin, $v_{i-1}$, and $v_{i+1}$.  Thus, replacing $v_i$ with $v'_i
$ if
necessary, we may assume that $L$ is convex.

Most of the remaining deformations will be described as pivots.  We
will fix an axis and rotate a vertex around a circle centered along
and perpendicular to the given axis.  If, for example, $|v_3-v_4|<k_3$,
we pivot the vertex $v_4$
around the axis through $0$ and $v_1$ until $|v_3-v_4|=k_3$.
It
follows from the choice of axes that the distances from $v_4$ to the origin
and $v_1$ are left unchanged, and it follows from the convexity of $L$
that $|w-v_4|$ increases for $w=v_2$, $v_3$, $v$, and $v'$.
    Similarly,
we may pivot vertices $v_i$ along the axis
through the origin and $v_{i+1}$ until $|v_i-v_{i-1}|=k_i$, for all $i$
.

Fix an axis through two opposite vertices (say $v_1$ and $v_3$)
and pivot another  vertex (say $v_2$) around the axis toward the
origin.  We wish to continue by picking different axes and
pivoting until $|v_i|=2$, for $i=1,2,3,4$.  However,
this process appears to break down in the event that
a vertex $v_i$ has distance $h_i$ from one of the enclosed
vertices and the pivot toward the origin moves $v_i$ closer
to that enclosed vector.  We must check that this situation can
be avoided.

Interchanging the roles of $(v_1,v_3)$ with
$(v_2,v_4)$ as necessary, we continue to pivot until
$|v_1|=|v_3|=2$ or $|v_2|=|v_4|=2$ (say the former).
We claim that either $v_2$ has distance greater than $h_2$ from
$v$ or that pivoting $v_2$ around
the axis through $v_1$ and $v_3$ moves $v_2$
away from $v$.  If not, we find that
$|v_1|=|v_3|=2$, $|v-v_2|=h_2$ and that $v$ lies in the
cone $C=C(v_2)$ with vertex $v_2$ spanned by the vectors
from $v_2$ to the origin, $v_1$, and $v_3$. (This
relies on the convexity of the region $L$.)

To complete the proof, we show that this figure made from $(0,v_1,v_2,v_3,v)$
cannot exist.  
Contract the edge $(v,v_2)$ as much as possible keeping 
the triangle $(0,v_1,v_3)$
fixed, subject to the constraints that $v\in C(v_2)$ and $|w-w'|\ge2$,
for $w=v,v_2$ and  $w'=0,v_1,v_3$.  This contraction gives
$|v-v_2|\le h_2<2\sqrt{2}$.

{\bf Case 1:} $v$ lies in the plane of $(0,v_2,v_3)$.  This gives an
impossibility: crossing edges $(v,v_2)$
$(0,v_3)$ of length
less than $2\sqrt{2}$.  Similarly, $v$ cannot lie in the plane of $(0,v_1,v_2)$.

{\bf Case 2:} $v$ lies in the interior of the cone $C(v_2)$. The contraction
gives
$|v-w'|=2$ for $w=v,v_2$ and $w'=0,v_1,v_3$.
The edge $(v,v_2)$ divides the
convex hull of  $(0,v_1,v_2,v_3,v)$
into three simplices.
Consider the dihedral angles of these simplices along
this edge. The dihedral angle of the simplex
$(v_1,v_2,v_3,v)$ is less than $\pi$.
The dihedral angles of the other two are less than
$\dih(S(2\sqrt{2},2,2,2,2,2))=\pi/2$.  Hence, the dihedral
angles along the diagonal cannot sum to $2\pi$
and the figure does not exist.

{\bf Case 3:} $v$ lies in the plane of $(v_2,v_1,v_3)$.
Let $r$ be the radius of the circle in this plane passing
through $v_1$ and $v_3$ obtained by intersecting the plane
with a sphere of radius $2$ at the origin.  
We have $2r\ge|v_1-v_3|>2\sqrt{2}$ because otherwise we have the impossible
situation of crossing edges $(v_1,v_3)$ and
 $(v_2,v)$ of length less than $2\sqrt{2}$.
Let $H$ be the perpendicular bisector of
the segment $(v_1,v_3)$.
By reflecting $v$ through $H$ if necessary
we may assume that $v$ and $v_2$ lie
on the same side of $H$, say the side of $v_3$.  Furthermore,
by contracting $(v,v_2)$, we may assume
without loss of generality 
that $|v_3-v|=|v_2-v_3|=2$.  
Let $f(r)=|v-v_2|$,
as a function of $r$.  $f$ is increasing
in $r$.  The inequalities
$2\sqrt{2}>h_2\ge|v-v_2|\ge f(\sqrt{2})=2\sqrt{2}$
give the desired contradiction.\qed

{\bf Proof (4.2):}
Assume for a contradiction that $v$ and $v'$ are vertices enclosed
by $L$.  Let the center of the Delaunay star be at the origin,
and let $v_1,\ldots, v_4$, indexed
consecutively,  be the four vertices of the Delaunay
star  that determine the extreme points
of  $L$.

We will describe a sequence of
deformations of the configuration 
(formed by the vertices $v_1,\ldots,v_4,v,v'$) that
transform the original configuration of vertices into particular
rigid arrangements below.  We will show that these rigid arrangements
cannot exist, and from this it will follow that the original configuration
does not exist either.
These deformations will preserve the constraints of the problem.
To be explicit, we assume that
that $2\le|w|\le 2.51$,
that $2\le |w-v_i|$ if $w\ne v_i$,
that $2\le |v-v'|$,
and that $|v_i-v_{i+1}|\le 2.51$,
for $i=1,2,3,4$ and $w=v_1,\ldots,v_4,v,v'$.
Here and
elsewhere we take our subscripts modulo 4, so that $v_1=v_5$, and
so forth.  The deformations will also keep $v$ and $v'$
in the cone at the origin that is determined by the
vertices $v_i$.  

We consider some deformations that increase $|v-v'|$.  
By Lemma 4.3, we may assume that
$|v_i|=2$,
$|v|=|v'|=|v_i-v_{i+1}|=2.51$, for $i=1,2,3,4$.
If, for some $i$,
we have $|v_i-v|>2$ and $|v_i-v'|>2$, then we fix $v_{i-1}$ and
$v_{i+2}$ and pivot $v_{i+1}$ around the axis through the origin and $v_{i+2}$
away from $v$ and $v'$.  The constraints $|v_{i+1}-v_i|=|v_i-v_{i-1}|=2.51$
will force us to drag $v_i$ to a new position on the sphere of radius 2.
By
making this pivot sufficiently small, we may assume that $|v_i-w|$ and
$|v_{i+1}-w|$, for $w=v,v'$, are greater than 2.

The vertices $v$ and $v'$ cannot both have distance 2 from both
$v_{i+2}$ and $v_{i-1}$, for then we would have $v=v'$.  So one of them,
say $v$, has distance exactly 2 from at most one of $v_{i+2}$ and $v_{i-1}$
(say $v_{i+2}$).  Thus, $v$ may be pivoted around the axis through the
origin and $v_{i+2}$
away from $v'$.  In this way, we increase $|v-v'|$ until
$|v_i-v|=2$ or $|v_i-v'|=2$, for $i=1,2,3,4$.

Suppose one of $v$, $v'$ (say $v$) has distance 2 
from $v_i$, $v_{i+1}$, and $v_{i+2}$.
The configuration is completely rigid. By symmetry, the vertices
$v_{i-1}$ and $v'$ must be the reflections of $v_{i+1}$
and $v$, respectively, through the plane through $0$, $v_i$,
and $v_{i+2}$.  In particular, $v'$ has distance 2 from
$v_i$, $v_{i-1}$, and $v_{i+2}$.
We pick coordinates and evaluate the
length $|v-v'|$.
We find that $|v-v'|\approx 1.746$, contrary to the hypothesis
that the centers of the spheres of our packing are separated by distances
of at least 2.
Thus, the
hypothesis that $v$ has distance two from three other vertices
is incorrect.

\gram|2|4.4|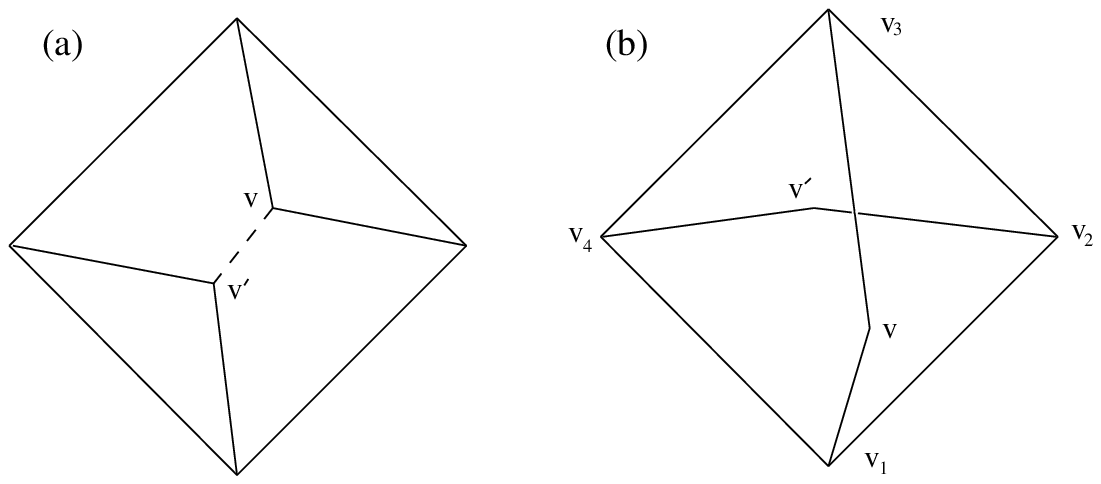|

We are left with one of the configurations of Diagram 4.4.  
An edge is drawn in the diagram,
when the distance between the two endpoints has the smallest possible value
(that is, $2.51$ for the four edges of the quadrilateral,
and $2$ for the remaining edges).  
Deform the figure of case (a) along the one remaining
degree of freedom until
$|v-v'|=2$.  In case (a),
referring to the notation established by Diagram 4.5,
we have a quadrilateral on the unit sphere
of edges $t_1 =  2\arcsin(2.51/4)\approx 1.357$ radians,
$t_2=\arccos(2.51/4)$ radians, and $t_3=2\arcsin(1/2.51)$ radians.
The form of this quadrilateral
is determined by the angle $\alpha$, and it is clear that
the angle $\beta(\alpha)$
is decreasing in $\alpha$.
 We have $$\theta=\dihmax =
\dih(S(2.51,2,2,2.51,2,2))\approx 1.874444.$$
The figure exists if and only if there exists $\alpha$ such that
$\beta(\alpha) + \beta(2\pi-\theta-\alpha) = 2\pi-\theta$.
By symmetry, we may assume that $\alpha\le\pi-\theta/2\approx 2.20437$.  
The condition $|v-v_1|\ge 2$ implies that
$$\alpha\ge \arccos\left (( \cos t_2 - \cos t_2 \cos t_3)/(\sin t_2 \sin t_3) \right)
        > 1.21.$$
However, by monotonicity,
$$\beta(\alpha) +\beta(2\pi-\theta-\alpha) < \beta(\alpha_{i}) +
        \beta(2\pi-\theta-\alpha_{i}-0.1)
        < 2\pi-\theta,$$
for $\alpha_{i}\le \alpha\le 0.1+\alpha_{i}$,
with $\alpha_i = 1.21+0.1\,i$, and $i=0,1,\ldots,9$, as a direct
calculation of the constants $\beta(\alpha_i) +\beta(2\pi-\theta-\alpha_i-0.1)$
will reveal.  (The largest constant,
which is about $2\pi-\theta-0.113$,  occurs for $i=0$.)
Hence, the figure of Diagram 4.4.a does not exist.

\gram|2|4.5|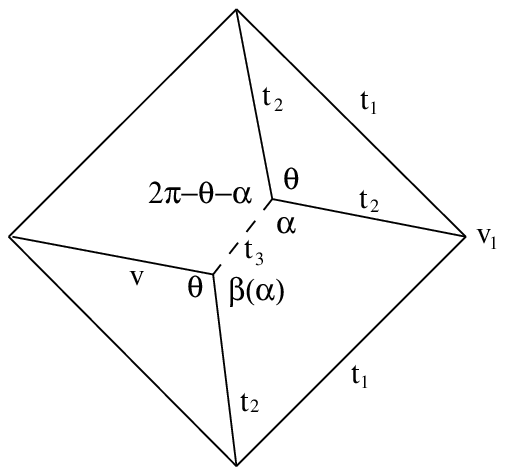|

To rule out Diagram 4.4.b, 
we reflect $v$, if necessary, to its image through the
plane $P$ through $(0,v_1,v_3)$, so that $P$
separates $v$ and $v'$.  The vertex $v$ can then be pivoted
away from $v'$ along the axis through $v_1$ and $v_3$.
This decreases $|v|$, but we may rescale so that $|v|=2.51$.
Eventually $|v-v_2|=2$ or $|v-v_4|=2$.  This is the
previously considered case in which $v$ has distance 2 from
three of the vertices $v_i$.
This completes the proof that the original arrangement of two enclosed
vertices does not exist. \qed

\bigskip
\centerline{\bf Section 5. Restrictions}
\bigskip

If a Delaunay star $D^*$ is composed entirely of tetrahedra, then
we obtain a triangulation of the unit sphere.  As explained in
Section 2, we wish to prove that no matter what the triangulation
is, we always obtain a score less than $8\,\pt$.
In this section we make a long list of properties that a configuration
must have if it is to have a score of $8\,\pt$ or more.  The next
section and the appendix show 
that only one triangulation satisfies all of these
properties.  Additional arguments will
show that this
triangulation scores less than $8\,\pt$.  This will complete the
proof of Theorem 1.

In this section, the term {\it vertices\/} refers to the vertices of the
triangulation of the unit sphere.
The edges of the triangulation give a planar graph.
We adopt the standard terminology of graph theory to describe the
triangulation.  We will speak of the degree of a vertex, adjacent
vertices, and so forth.
The $n$ triangles around a vertex will be
referred to as an $n$-gon.  We will also refer to the corresponding
$n$ tetrahedra that give the triangles of the polygon.
We say that a triangulation contains
a {\it pattern\/} $(a_1,\ldots,a_n)$, for $a_i\in {\Bbb N}$, 
if there are distinct
vertices $v_i$ of degrees $a_i$ that are  pairwise nonadjacent, 
for $i=1,\ldots,n$. Let $N$ be the number of vertices in the
triangulation, and let $N_i$ be the number of vertices of degree $i$.
We have $N = \sum N_i$.

In this section, we will start to use various inequalities
from Section 9 related to the score.  Since the score $\score(S)$
may be either $\Gamma(S)$ or $\vor(S)$ depending on the
circumradius of $S$, there are two cases to consider for
every inequality.  In general, the inequalities for $\Gamma(S)$
are more difficult to establish.  In the following sections,
we will only cite the inequalities pertaining to $\Gamma(S)$.
Section 9 shows how all the same inequalities hold
for $\vor(S)$.

{\bf Proposition 5.1.} Consider a Delaunay star $D^*$ that is
composed entirely of quasi-regular tetrahedra.  
Suppose that $\score(D^*)\ge 8\,\pt$.  Then
the following restrictions hold
on the triangulation of the unit sphere given by
the standard decomposition.

{
\def\ha{\hangindent=20pt\hangafter=1\relax}

\parskip=0.3\parskip
\hbox{}

\ha
1. $13\le N\le 15$.

\ha
2.  $N = N_4+N_5+N_6$.

\ha
3.  A region bounded by three edges is either a single triangle or the
	complement of a single triangle.

\ha
4. Two degree four vertices cannot be adjacent.

\ha
5. $N_4\le 2$.

\ha
6. Patterns $(6,6,6)$ and $(6,6,4)$ do not exist.

\ha
7. The pattern $(6,6)$ or $(6,4,4)$ implies that $N\ge 14$.

\ha
8. If there are two adjacent degree six vertices, and a third degree six
	vertex adjacent to neither of the first two, then $N=15$ and
	all other vertices are adjacent to at least one of these three.

\ha
9. The triangulation is made of geodesic arcs on the sphere whose radian
lengths are between $0.8$ and $1.36$.

}

\bigskip

{\bf Lemma 5.2.}  Consider a vertex of degree $n$,
for some $4\le n\le 7$.
Let $S_1,\ldots, S_n$ be the tetrahedra
that give the $n$ triangles.   Then $\sum_{i=1}^n \score(S_i)$ is
less than $z_n$, where $z_4 = 0.33\,\pt$, $z_5 = 4.52\,\pt$,
$z_6 = -1.52\,\pt$, and $z_7 = -8.9\,\pt$.
Suppose that $n\ge 6$.  Let $S_1,\ldots,S_4$ be any four of
the $n$ tetrahedra around the vertex.  Then $\sum_{i=1}^4\score(S_i)
        < 1.5\,\pt$.

\bigskip
{\bf Proof (5.2):}
When $n=4$, this is Lemma 9.5.  When $n=5$, this is
Lemma 9.6.  When $n=6$ or $n=7$, we have by Calculation 9.4
$$
\align
\sum_{i=1}^n \score(S_i) &< \sum_{i=1}^n
        (0.378979\dih(S_i)  -0.410894) \\
	&= 2\pi(0.378979)-0.410894n.
\endalign$$
The right-hand side evaluates to about $-1.520014\,\pt$ and $-8.940403\,\pt$,
respectively,
when $n=6$ and $n=7$.

Assume that $n\ge 6$, and select any four $S_1,\ldots,S_4$ of
the $n$ tetrahedra around the vertex.  Let $S_5$ and $S_6$ be
two other tetrahedra around the vertex.
The dihedral angles of $S_5$ and $S_6$ are are at least $\dihmin$,
by Calculation 9.3.  Each of the four triangles (associated
with $S_1,\ldots,S_4$) must then, on average, have an
angle at most $(2\pi - 2\dihmin)/4$ at $v$.  By Calculation 9.4,
$$
\align
\sum_{i=1}^4 \score(S_i) &< \sum_{i=1}^4 (0.378979 \dih(S_i)  -0.410894)\\
&\le (2\pi - 2\dihmin) 0.378979 +4 (-0.410894) < 1.5\,\pt.
\endalign$$\qed

\smallskip
{\bf Proof (5.1):}  
Let $t = 2(N-2)$ be the number of triangles, an even number.
By Euler's theorem
on polyhedra,
$$3 N_3 + 2 N_4 + N_5 + 0 N_6 - N_7 \cdots = 12.$$
Let $S_1,\ldots,S_t$ denote the tetrahedra of $D^*$.
Let $\score_i=\score(S_i)$ be the corresponding score.
Let $\sol_i$ denote the solid angle
cut out by $S_i$ at the origin.  We have $\sum \sol_i = 4\pi$.
Often, without warning, we will rearrange the indices $i$ so
that the tetrahedra that give the triangles around a given vertex
are numbered consecutively.  When a vertex $v$ of the triangulation
has been fixed, we will let $\alpha_i$ denote the angles of the
triangles at $v$, so that $\sum\alpha_i = 2\pi$. The angle
$\alpha_i$ of a triangle is equal to the
corresponding dihedral angle of the simplex
$S_i$.
We abbreviate certain sums over $n$ elements to $\sum_{(n)}$, when
the context makes the indexing set clear.
Throughout the argument, we will use Calculation 9.1, which asserts
that $\score_i\le 1\,\pt$, for all $i$. The proofs will show that
if a triangulation fails to have any of the properties 1--9,
then the total score
must be less than $8\,\pt$.

(Proof of 5.1.1):
Assume that $t\ge 28$.  By Calculation 9.9,
$$
\align
\sum_{(t)}\score_i &< \sum_{(t)} (0.446634 \sol_i  -0.190249)
        \le {4\pi} (0.446634) + t (-0.190249)\\
&\le {4\pi} (0.446634) + 28 (-0.190249) \le 8\,\pt.
\endalign$$

Suppose that $t\le 18$.  By Calculation 9.8,
$$\align
\sum_{(t)} \score_i &< \sum_{(t)}
(-0.37642101 \sol_i + 0.287389) = {4\pi} (-0.37642101) + t (0.287389)\\
&\le {4\pi} (-0.37642101) + 18(0.287389) \le 8\,\pt.
\endalign
$$
This proves that $12\le N\le 15$.  The case $N=12$ will be excluded after
5.1.2.

(Proof of 5.1.2):  By Lemma 8.3.2 
 and Calculation 9.3, each angle $\alpha$
of each triangle satisfies $\dihmin\le \alpha\le \dihmax$.  In
particular, $\alpha>\pi/4$, so that each vertex has degree
less than eight,
 and $\alpha<2\pi/3$, so that each vertex has degree greater
than three.

Assume for a contradiction that there
is a vertex of degree seven.  Consider the seven
tetrahedra around the given vertex.  By Lemma 5.2, the tetrahedra
satisfy
$$\sum_{(7)} \score_i <  -8.9\,\pt.$$
Suppose $t\le 24$.
For each vertex $v$,
set
$\zeta_v = \sum_{ v }\score_i$, the sum running over the
tetrahedra around $v$.
Clearly, $\sum_{(t)}\score_i = (\sum_v \zeta_v)/3$.
Pick a vertex $v$ that is not a vertex
of any of the seven triangles of the heptagon.
By Lemma 5.2, we see that $\zeta_v < 0.33\,\pt$
if $v$ has degree four,
$\zeta_v < 4.52\,\pt$ if $v$ has degree five, 
$\zeta_v< -1.52\,\pt$ if $v$ has degree six,  and
$\zeta_v < -8.9\,\pt$ if $v$ has degree seven.
In particular, if $v$ has degree $n$, then $\zeta_v$ falls
short of $n$ points by at least $0.48\,\pt$.
Thus,
$$\sum_{(t)} \score_i < \sum_{(7)}\score_i + \sum_{(t-7)}\,\pt - 0.48\,\pt
\le (-8.9 + (24-7) -0.48)\,\pt < 8\,\pt.$$

Finally, we assume that $t=26$ and $N=15$.  If $N_6=0$ and $N_7=1$,
then Euler's theorem gives the incompatible conditions
$2N_4 + N_5 -1=12$ and $N_4+N_5+1=15$. Thus, $N_6>0$ or $N_7>1$.
This gives a second $k$-gon $(k=6 \hbox{ or } 7)$ around a vertex $v$.
This second polygon
shares at most two triangles with the original heptagon.  This leaves
at least four triangles of a second polygon exterior to the first.
Thus, by Lemma 5.2,
$$\sum_{(26)} \score_i
        \le \sum_{(7)}\score_i +\sum_{(4)}\score_i +
        \sum_{(t-11)}\,\pt < -8.9\,\pt + 1.5\,\pt + (26-11) \,\pt < 8\,\pt.$$

(Proof of 5.1.1, cont.):  Assume $N=12$.
We define three classes
of quasi-regular tetrahedra.  
In the first class, all the edges have lengths between $2$ and $2.1$.
In the second, the fourth, fifth, and sixth edges
have lengths greater than $2.1$.  The third class is everything else.

Set $\epsilon = 0.001$, $a=-0.419351$, $b= 0.2856354$.
The following are established
by Calculations 9.10 -- 9.12.
$$\align
	\score(S) &\le a \sol(S) + b + \epsilon,
	\quad \hbox{for $S$ in the first class,} \\
	\score(S) &\le a \sol(S) + b,
	\quad \hbox{for $S$ in the second class,} \\
	\score(S) &\le a \sol(S) + b -5 \epsilon,
	\quad \hbox{for $S$ in the third class,} \\
\endalign
$$
Consider a vertex $v$ of degree $n=4$, $5$, or $6$
and the surrounding
 tetrahedra
$S_1,\ldots, S_n$. 
We claim that 
$$\sum_{i=1}^n \score(S_i) \le \sum_{i=1}^n (a\sol(S_i) +b).\tag5.1.1.1$$ 
This follows directly from the stated inequalities
if none of these tetrahedra
are in the first class.  It
is also obvious if at least one of these tetrahedra is in the third class,
because then the inequality is violated by at most 
$-5\epsilon + (n-1)\epsilon \le0$.
So assume that all of the tetrahedra are in the first two classes with
at least one in the first class.  By the restrictions on the lengths of
the edges, a tetrahedron in the first class cannot be adjacent to
one in the second class.  We conclude that the tetrahedra are all in the
first class.

By Lemma 5.2, $z_n \le n(0.904)\pt$.
If $\sum_{(n)} \sol(S_i) \le n(0.56176)$, then
$$\sum_{(n)}\sigma(S_i) \le n (0.904)\,\pt \le \sum_{(n)} (a (0.56176) + b) \le
	\sum_{(n)} (a\sol(S_i) + b).$$
So we may assume that $\sum_{(n)}\sol(S_i) \ge n(0.56176)$.
By Calculation 9.13, 
$$\align
	\sum_{(n)} \score(S_i) &\le \sum_{(n)}(-0.65557 \sol(S_i) + 0.418) \\
	\le
	\sum_{(n)}& (a\sol(S_i) + b) + \sum_{(n)} (0.132365 - 0.236219 (0.56176)) 
	< \sum_{(n)} (a\sol(S_i) + b).
\endalign
$$
This establishes inequality 5.1.1.1.  By averaging over every vertex,
we see that the average of the scores $\score(S)$ is less than
the average of $a\sol(S) + b$.  So
$$\score(D^*)= \sum_{(t)} \sigma(S) \le
	\sum_{(t)} (a\sol(S) + b) = 4\pi a + 20 b \approx 7.99998\,\pt < 8\,\pt.$$

(Proof of 5.1.3):  This is Lemma 3.7.  (If we had
not introduced quasi-regular tetrahedra, then this result would
no longer hold.)

(Proof of 5.1.4):  If two degree four vertices
are adjacent,
we then have the
arrangement of Diagram 5.3.  This is a quadrilateral enclosing
two vertices, contrary to Proposition 4.2.

\gram|2|5.3|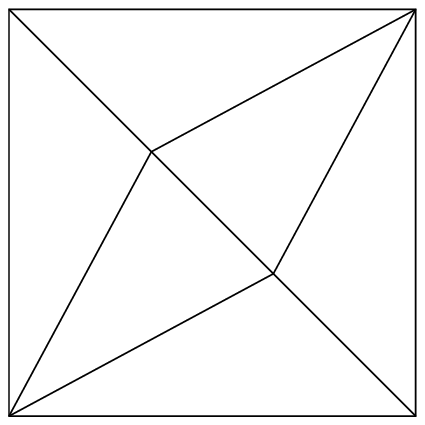|  

(Proof of 5.1.5):  As in the proof of 5.1.2, for each vertex $v$, set
$\zeta_v = \sum_{ v }\score_i$, the sum running over the
tetrahedra around the vertex $v$.  We use the estimates of $\zeta_v$
that appear in Lemma 5.2.

We have found that $N_i=0$, if $i\ne 4,5,6$.  By Euler's theorem,
$N_5 = 12-2N_4$ and $N_6 = N_4 + N-12$.  Assume that $N\ge 13$ and
that $N_4\ge3$.  Then
$$\align
\sum_{(t)}\score_i &= {1\over 3}\sum_{(N)}\zeta_v
        < {1\over3}(0.33 N_4 + 4.52 (12-2N_4) -1.52 (N_4+N-12))\,\pt
\\
&\le {1\over 3}(0.33 (3) + 4.52 (12-6) -1.52 (3+13-12))\,\pt < 8\,\pt.
\endalign
$$

(Proof of 5.1.6):  Suppose that we have the pattern $(6,6,6)$.  Then
reordering indices according to the polygons in the pattern,
we have by the estimates of Lemma 5.2
$$\align
\sum_{(t)}\score_i
&\le \sum_{(6)} \score_i + \sum_{(6)}\score_i + \sum_{(6)}\score_i
                + \sum_{(t-18)} \score_i \\
&< -1.52\,\pt - 1.52 \,\pt - 1.52 \,\pt + (26-18)\,\pt  < 3.5\,\pt.
\endalign
$$
Similarly, if we have the pattern $(6,6,4)$, then we find
$$\sum_{(t)}\score_i \le -1.52\,\pt -1.52\,\pt + 
	0.33\,\pt + (26-16)\,\pt < 8\,\pt.$$

(Proof of 5.1.7):  We use the same method as in the proof of
5.1.6.  If we have the pattern $(6,6)$, 
and if $t\le 22$, then
$$\sum_{(t)}
\score_i \le \sum_{(6)}\score_i + \sum_{(6)}\score_i + \sum_{(t-12)}\score_i
\le -1.52\,\pt - 1.52\,\pt + (22-12)\,\pt < 8\,\pt.$$
Similarly, if we have the pattern $(6,4,4)$,
then $\sum_{(t)}\score_i$ is less than
$-1.52\,\pt+0.33\,\pt+0.33\,\pt+8\,\pt < 8\,\pt$.

(Proof of 5.1.8):   Let the two adjacent degree six
vertices be $v_1$ and $v_2$.  Let the third be $v_3$.
The six
triangles in the hexagon around $v_3$ give less than $-1.52\,\pt$.
The ten triangles in the hexagons around $v_1$ and $v_2$ give at most
$$\sum_{(4)}\score_i + 
	\sum_{(6)} \score_i < 1.5\,\pt - 1.52\,\pt <0\,\pt,$$
by the argument described in the case $t=26$ of 5.1.2  (see Lemma 5.2).

Suppose that $t\le 24$.
There remain at most $24-16=8$ triangles,
and they give a combined score of at most $8\,\pt$.  The total score
is then less than $(-1.52+8)\,\pt$, as desired.

Now assume that $t=26$.
 Suppose there is a vertex $v$ that is not
adjacent to
any of $v_1$, $v_2$, or $v_3$. As in the proof of 5.1.8,
the ten triangles in the two overlapping hexagons
give less than $0\,\pt$.  The other hexagon gives less than $-1.52\,\pt$.
By Lemma 5.2, the $n$ triangles around $v$ fall short of $n$ points
by at least $(5-4.52)\,\pt=0.48\,\pt$.
Each of the remaining
triangles gives at most $1\,\pt$.  The score is then
less than $(0-1.52+(26-16)-0.48)\,\pt = 8\,\pt$, as desired.

(Proof of 5.1.9):  
This follows directly from the construction
of the triangulation and the close-neighbor
restrictions  
on the lengths of the
edges of a quasi-regular tetrahedron.
The lengths are between $0.8<2\arcsin(1/2.51)$ and $2\arcsin(2.51/4)<1.36$.
This completes the proof of
the proposition.\qed

\bigskip
\centerline{\bf Section 6. Combinatorics}
\bigskip
{\bf Theorem 6.1.}  Suppose that a triangulation satisfies Proposition 4.2
and Properties
1--9 of Proposition 5.1.
Then it must the 
triangulation of Diagram 6.2 with 14 vertices and 24 triangles.
\smallskip

\gram|1.5|6.2|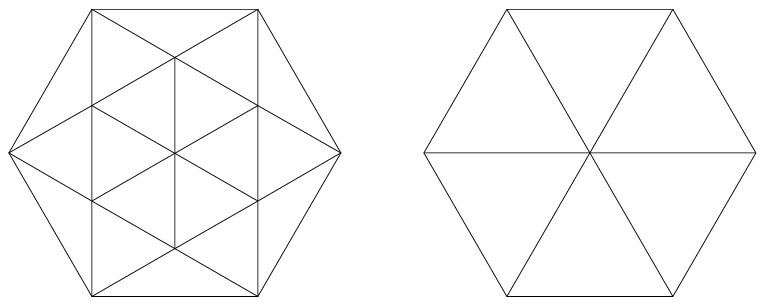|

{\bf Proof:}
Fix a polygon centered at a vertex $u_0$, such as the
hexagon  in Diagram 6.3.
The six vertices
$v_1,\ldots,v_6$ of the hexagon are distinct, for otherwise two
distinct geodesic arcs on the sphere would run between $u_0$ and
$v_i$ for some $i$.  This is impossible, because $u_0$ and $v_i$
are not antipodal by Property 5.1.9.  Similarly, the vertices of every
other polygon of the triangulation are distinct.

\smallskip
\gram|1.5|6.3|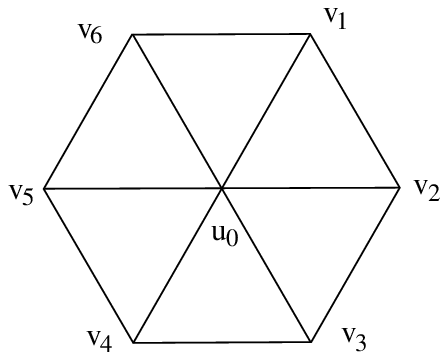|
\smallskip

We then extend the polygon to a second layer of triangles.
Each of the vertices $v_i$ has degree four, five, or
six,
and two degree four vertices cannot be adjacent.
The new vertices are denoted $w_1,\ldots,w_k$.
One example is shown in Diagram 6.4.

\gram|2|6.4|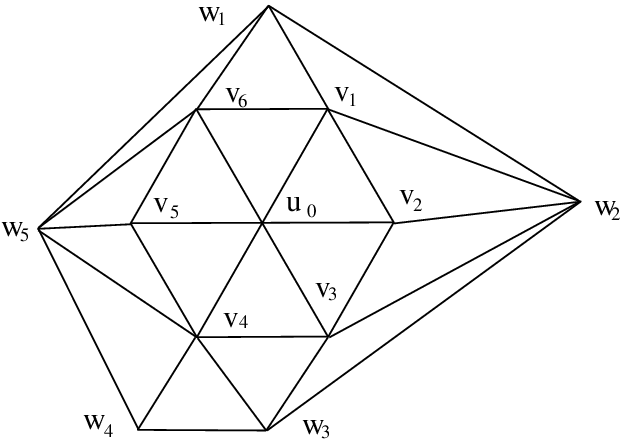| 

There can be no identification of a vertex $v_i$ with a vertex $w_j$,
for otherwise there is a triangle (say with vertices
$w_j,v_k,u_0,v_i=w_j$)
that is subtriangulated, contrary to Property 5.1.3.
Similarly, there
is no identification of two vertices $w_i$ and $w_j$, for otherwise
it can be checked that there is
a quadrilateral (with vertices $w_i,v_k,u_0,v_\ell,w_j=w_i$)
that encloses more than one vertex, which is impossible by Proposition
4.2.  A purely combinatorial problem remains.  It is solved
in the appendix. \qed

\bigskip
{\bf Proposition 6.5.}  The triangulation of Theorem 6.1 scores
less than $8\,\pt$.

{\bf Proof:}
Our initial bound on the score comes by viewing the
triangulation as made up of two hexagons and twelve additional
triangles.  By Lemma 5.2 and Calculation 9.1,
$$\sum_{(24)} \score_i \le \sum_{(6)}\score_i +
\sum_{(6)}\score_i + \sum_{(12)} \,\pt 
	< -1.52\,\pt -1.52\,\pt + 12\,\pt = 8.96 \,\pt.
\tag6.6$$
A refinement is required in order to lower the upper bound to $8\,\pt$.

If a vertex $v$ has height $|v|\ge2.2$, 
then by Calculation 9.2, we find
$\score_i < 0.5\,\pt$ for the tetrahedra at $v$.  Thus, in the Inequality
6.6, if the vertex $v$ of some pentagon has height 
$|v|\ge 2.2$, then the
term $\sum_{(12)}\,\pt$\thinspace\
may be replaced by $\sum_{(9)}\,\pt + \sum_{(3)} 0.5\,\pt$
\thinspace\ (there are
three triangles in the pentagon that do not belong to the hexagon),
and the upper bound on
the score falls to $7.46\,\pt$.

If a vertex of degree six has height $|v|\ge 2.05$, 
then we claim that $\sum_{(6)}\score_i < -3.04\,\pt$.
In fact,  by Calculation 9.7, the hexagon gives
$$\sum_{(6)} \score_i < \sum_{(6)} (0.389195\dih(S_i) -0.435643) =
                2\pi (0.389195) - 6(0.435643) < - 3.04\,\pt.$$
Estimate 6.6 is improved to
$$\sum_{(24)} \score_i < -3.04\,\pt -1.52\,\pt + 12\,\pt  = 7.44\,\pt.$$

If the twelve tetrahedra have combined solid angle less than $6.48$, then
by Calculation 9.9,
$$\align
\sum_{(12)}\score_i &< \sum_{(12)} (0.446634\sol(S_i)  -0.190249)\\
        &\le 6.48(0.446634)+12(-0.190249)
< 11.039\,\pt.
\endalign$$
Then Estimate 6.6 is improved to the bound
$-1.52\,\pt-1.52\,\pt + 11.039\,\pt = 7.999\,\pt$.

Now assume, on the other hand,
 that the combined solid angle of the two hexagons is at
most $4\pi- 6.48$.  Set $K = (4\pi-6.48)/12$ and define
$\score'(S):= \score(S) + (K -\sol(S))/3$.
The solid angle of one of the two
hexagons is at most $6K$.
For that hexagon, we have
$$\sum_{(6)} \score'(S_i) = \sum_{(6)} \score_i +
        (6K - \sum_{(6)}\sol(S_i))/3
        \ge \sum_{(6)}\score_i.$$
By our previous estimates, 
we now assume without loss of generality that the heights $|v|$ of
the vertices of triangles
in the hexagon are at most $2.05$, $2.2$, and $2.2$, the bound
of $2.05$ occurring at the center of the hexagon.
By Calculation 9.15,
$$\sum_{(6)} \score'(S_i) \le \sum_{(6)} 0.564978\dih_i-0.614725
        = 2\pi(0.564978)+6(-0.614725) < -2.5\,\pt.$$
Estimate 6.6 becomes
$$\sum_{(24)}\score_i \le \sum_{(6)}\score'(S_i) + \sum_{(6)}\score_i +
        \sum_{(12)}\score_i < (-2.5-1.52+12)\,\pt = 7.98\,\pt.$$
This completes the proof of the proposition.\qed

\bigskip
\centerline{\bf Section 7. The Method of Subdivision}
\bigskip
The rest of this paper is devoted to the verification of the
inequalities that have been used in Sections 5 and 6.
In this section we describe the method used to obtain several of
our bounds.  We call it the method of subdivision.
Let $p(x) = \sum_I c_I x^I$ be a polynomial with real coefficients
$c_I$, where
$I = (i_1,\ldots,i_n)$, $x=(x_1,\ldots,x_n)\in {\Bbb R}^n$, and
$x^I = x_1^{i_1}\cdots x_n^{i_n}$.
It is clear that if $C$ is the product of intervals $[a_1,b_1]\times
\cdots \times [a_n,b_n] \subset {\Bbb R}^n$ in the positive orthant
($a_i>0$, for $i=1,\ldots,n$), then
$$\forall x\in C,\qquad
p_{\min}(C) \le p(x) \le p_{\max}(C),$$
where $$p_{\min}(C) = \sum_{c_I>0} c_I a^I + \sum_{c_I<0} c_I b^I
\quad\hbox{and}\quad
p_{\max}(C) = \sum_{c_I>0} c_I b^I + \sum_{c_I<0} c_I a^I.$$
Another bound comes from the Taylor polynomial $p(x) = \sum d_I
(x-a)^I$ at $a$:
$$d_0 + \sum_{d_I<0,I\ne0}
   d_I(b-a)^I \le p(x) \le d_0 +\sum_{d_I>0,I\ne0}
    d_I(b-a)^I.\tag 7.1$$
If $r(x)=p(x)/q(x)$ is a rational function, and
if $q_{\min}(C)>0$, then
$$\forall x\in C,\qquad
{p_{\min}(C)\over q_1(C,p)} \le r(x) \le {p_{\max}(C)\over q_{2}(C,p)},$$
where
$q_1(C,p)$ (resp. $q_2(C,p)$)
is defined as $q_{\max}(C)$ whenever $p_{\min}(C) \ge 0$ (resp. $p_{\max}(C)<0$)
and as $q_{\min}(C)$ otherwise.

Let us define a {\it cell\/} to be a product of intervals
in the positive orthant of ${\Bbb R}^n$.
By covering a region with a sufficiently fine collection of cells,
 various inequalities of rational functions are easily
established.
To prove an inequality of rational functions with
positive denominators (say $r_1(x) < r_2(x)$, for all $x\in C$), we
cover $C$ with a finite number of cells and compare the upper bound
of $r_1(x)$ with the lower bound of $r_2(x)$ on each cell.  If it turns
out that some of the cells give too coarse a bound, then we
subdivide each of the delinquent cells into a number of smaller
cells and repeat the process.  If at some stage we succeed in
covering the original region $C$ with cells on which the upper bound
of $r_1(x)$ is less than the lower bound of $r_2(x)$, the inequality
is established.

A refinement of this approach applies the method to the partial
derivatives.
  If, for instance, we establish by the method of
subdivision that for some $i$
$${\partial p\over \partial x_i}(x) \ge 0,\qquad \forall x\in C,$$
then we may compute an upper bound of $p$ by applying the method
of subdivision to the polynomial obtained from $p$ by the specialization
$x_i=b_i$, where $b_i$ is the upper bound of $x_i$ on $C$.
Thus, we obtain an upper bound on a polynomial by fixing all the
variables that are known to have partial derivatives
  of fixed sign and then
applying the method subdivision to the resulting polynomial.
Similar considerations apply to lower bounds and to
rational functions with
positive denominators.

It is a fortunate circumstance that many of the polynomials we
encounter in sphere packings are quadratic in each variable with
negative leading coefficient
($u$, $\rho$, $\Delta$, $\chi$ in the  next section).  
In this case, the lower bound is attained at a corner of the
cell.  Of course, the maximum of a quadratic function with
negative leading coefficient is also elementary:
$-\alpha(x-x_0)^2+\beta\le\beta$, if $\alpha\ge0$.

\bigskip
\centerline{\bf Section 8. Explicit Formulas for Compression, Volume,
        and Angle}
\smallskip
Many of the formulas in this section are classical.  One can
typically find them in nineteenth century primers on solid
geometry.  The formula for solid angles, for example, is due
to Euler and Lagrange.  For anyone equipped with symbolic
algebra software, the verifications are elementary, so we omit
many of the details.
All formulas in this
section will be valid for Delaunay simplices and for quasi-regular
tetrahedra, unless otherwise noted.

{\bf 8.1. The Volume of a Simplex}

As in the previous section, a {\it cell\/} in ${\Bbb R}^n$ is a 
product of intervals in ${\Bbb R}^n$.
We define a function $\Delta:[4,16]^6 
	\subset {\Bbb R}^6 \to {\Bbb R}$ by
$$\align
\Delta(x_1,\ldots,x_6) = &x_1 x_4(-x_1+x_2+x_3-x_4+x_5+x_6)\\
                        &+x_2 x_5(x_1-x_2+x_3+x_4-x_5+x_6)\\
                        &+ x_3 x_6 (x_1+x_2 -x_3 +x_4 +x_5 -x_6)\\
                        &-x_2 x_3x_4 - x_1 x_3 x_5 - x_1 x_2 x_6 - x_4 x_5 x_6.
\tag8.1.1
\endalign
$$
We set $y_i = \sqrt{x_i}$, for $i=1,\ldots, 6$. This relationship
between $x_i$ and $y_i$ remains in force to the end of the paper.
Index the edges of a simplex as in Diagram 8.1.2.

\gram|2|8.1.2|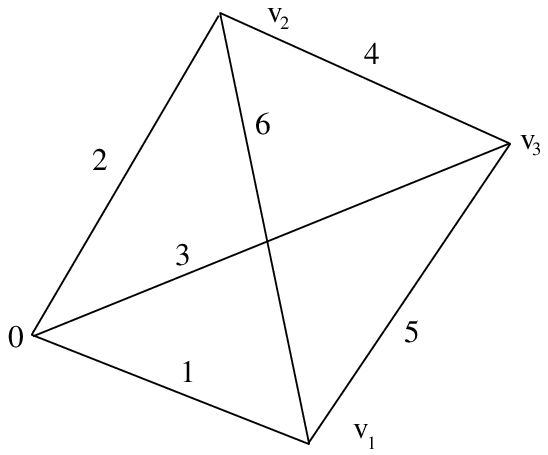|

We also define $u:[4,16]^3\to {\Bbb R}$ by
$$\align
u(x_1,x_2,x_6) &= (y_1+y_2+y_6)(y_1+y_2-y_6)(y_1-y_2+y_6)(-y_1+y_2+y_6)\\
&=-x_1^2 - x_2^2 - x_6^2 + 2 x_1 x_6 + 2 x_1 x_2 + 2 x_2 x_6.
\tag8.1.3\endalign
$$

{\bf Lemma 8.1.4.}  
There exists a simplex of positive volume with edges of length
$y_1,\ldots,y_6$ if and only if $\Delta(x_1,\ldots,x_6)>0$.  If these conditions
hold, then the simplex has volume $\Delta(x_1,\ldots,x_6)^{1\over2}/12$.

{\bf Proof:}
The function $u=u(x_1,x_2,x_6)$
is quadratic in each variable, with negative leading coefficient,
so the minimum of $u$, which is 0,
 is attained at a vertex of the cube $[4,16]^3$.  At the vertices
where the minimum is attained, $\Delta(x_1,\ldots,x_6)\le0$.

Assume $\Delta>0$.  Then $u>0$, and 
a simplex exists with vertices 
$0$, $X=(y_1,0,0)$, $Y={1\over2}(*,u^{1\over2}/y_1,0)$,
$Z=(*,*,{(\Delta/u)^{1\over2}})$.  
Conversely, if the simplex exists, then it must be of the given form, up to
an orthogonal transformation,
so $u>0$ (otherwise $O$, $X$, and $Y$ are colinear) and $\Delta> 0$.
The volume is $|\det(X,Y,Z)|/6  = \Delta^{1\over2}/12$.\qed

Let $C$ be a cell contained in $[4,16]^6$.  The minimum of $\Delta$
on $C$ is attained at a vertex of $C$: this is clear
because $\Delta$ is quadratic in each variable with negative
leading coefficient. To obtain an upper bound on $\Delta$ on a
cell, we may use the method of Section 7.  
The restriction of $\Delta$ to the zero set of $\partial\Delta/
\partial x_1$ is
$$ {u(x_4,x_5,x_6)\cdot u(x_2,x_3,x_4)
                        \over 4 x_4}.
$$
This is an upper bound on $\Delta(x_1,\ldots,x_6)$.

\bigskip
\line{\bf 8.2. The Circumradius of a Simplex\hfill}
\bigskip

The circumradius of a face with edges $y_4$, $y_5$, and $y_6$ is
$$\eta(y_4,y_5,y_6) = {y_4y_5y_6\over u(y_4^2,y_5^2,y_6^2)^{1/2}}.$$

Define $\rho:[4,16]^6 \to {\Bbb R}$ by
$$\rho(x_1,\ldots,x_6) = -x_1^2 x_4^2 - x_2^2 x_5^2 - x_3^2 x_6^2 +
        2 x_1 x_2 x_4 x_5 + 2 x_1 x_3 x_4 x_6 + 2 x_2 x_3 x_5 x_6.\tag8.2.1$$

{\bf Lemma 8.2.2.}  Suppose $\Delta>0$; then $\rho>0$.

{\bf Proof:}  If $\Delta>0$, then the simplex of Lemma 8.1.4 exists.
Let the coordinates of the circumcenter 
of the simplex be $(x,y,z)$.
Direct calculation shows that
$0< x^2+y^2+z^2 = \rho/(4\Delta)$.\qed

{\bf Corollary 8.2.3.}  Let $S$ be a Delaunay simplex
of positive volume with edges of lengths $y_i=\sqrt{x_i}$.
The  circumradius is
${1\over 2}(\rho/\Delta)^{1\over2}$.  

Set
$$
\align
\chi(x_1,x_2,x_3,x_4,x_5,x_6)
&= x_1 x_4 x_5 + x_1 x_6 x_4 + x_2 x_6 x_5 + x_2 x_4 x_5 + x_5 x_3 x_6 \\
        &\quad+ x_3 x_4 x_6 - 2 x_5 x_6 x_4 - x_1 x_4^2 - x_2 x_5^2 - x_3 x_6^2.
\endalign
$$
We have $${\rho\over 4\Delta} - \eta(y_4,y_5,y_6)^2 = 
{\chi(x_1,\ldots,x_6)^2\over 4u(x_4,x_5,x_6)\Delta}.$$
The vanishing of $\chi$ is the condition for the
circumcenter of the simplex to lie
in the plane through the face $T$
bounded by the fourth,
fifth, and sixth edges.  $\chi$ is positive if the circumcenter of
$S$ and the vertex of $S$ opposite $T$ lie on the same side of the
plane through $T$ and negative if they lie on opposite sides
of the plane.  We will say that $T$ has {\it positive orientation\/}
when $\chi>0$.

{\bf 8.2.4.}\ The function
 $\rho$ is quadratic in each of the variables $x_1,\ldots,x_6$ with
negative leading coefficient. Thus,
the minimum of $\rho$ is attained at a vertex of a given cell $C$.
The derivative is
$\Delta^2 \partial(\rho/\Delta)/\partial x_1 = 
\chi(x_5,x_6,x_1,x_2,x_3,x_4)\chi(x_1,x_2,x_3,x_4,x_5,x_6)$.
This leads to bounds on the circumradius by the method of
Section 7.

The circumradius of a quasi-regular tetrahedron is increasing
in $x_1,\ldots,x_6$ if the orientation of each face is positive.
When a face fails to have positive orientation, it satisfies
the constraints of Section 3 (for example $y_1,y_2,y_3\in [2,2.15]$,
$y_4,y_5,y_6\in [2.3,2.51]$, $\eta(y_4,y_5,y_6)\ge\sqrt{2}$, if
it is the face opposite the origin).
This allows us to determine bounds on the circumradius for
most of the quasi-regular tetrahedra
 we encounter in this paper by inspection.
For example, in Section 9 we will study the simplices constrained
by $y_i\in [2,2.1]$.  An upper bound on the circumradius is
$\rad(S(2.1,2.1,2.1,2.1,2.1,2.1))$.

\noindent
{\bf 8.2.5.}\quad Proof of Lemma 3.4.  In the notation of Section 3.4, 
let $T$ be the given face, and let $v_0$ be a vertex satisfying
the conditions of the lemma.
Let $S$ be the simplex at the origin $v_0=0$, whose first, second,
and third edges abut at $v_i$, for $i=1,2,3$.
Suppose for a contradiction that $|v_0-v_3|^2=x_3\ge 2.15^2$.
In light of the results of this section, a contradiction
follows if 
$\chi(x_1,\ldots,x_6)>0$, for
$x_i\ge 4$ for $i=1,2$, $x_3\ge 2.15^2$, $2.3^2\le x_i\le 2.51^2$
for $i=4,5,6$.
(The lower bound of 2.3 comes from Remarks
3.2 and 3.3: the circumradius of $T$ is at least $\sqrt{2}$).
Since $T$ is acute,
$${\partial\chi\over\partial x_1} = x_4 (-x_4+x_5+x_6)>0.$$
Similarly, $\chi$ is increasing in $x_2$ and $x_3$.  Thus, to
minimize $\chi$, we take $x_1=x_2=4$, $x_3=2.15^2$.
But $\chi$ is quadratic in each variable $x_4$, $x_5$, and $x_6$
with negative leading coefficient, so the minimum occurs at a corner
point of $[2.3^2,2.51^2]^3$.
We find $\chi\ge\chi(2^2,2^2,2.15^2,2.51^2,2.51^2,2.51^2)
\approx 0.885>0$.\qed

\medskip
For the set of Delaunay
stars to be a
compact topological space, simplices of zero volume
must be included.  
Since the circumradius remains bounded, the degenerate
Delaunay simplices of zero
volume are planar quadrilaterals that possess a circumscribing
circle.  

\bigskip

\line{\bf 8.3. Dihedral Angles\hfill}

Let $\dih(S)$ be the dihedral
angle of a simplex $S$
along the first edge.
It is the (interior) angle formed by
faces with edges $(2,1,6)$ and $(1,3,5)$.

{\bf Lemma 8.3.1.}
$$\cos\dih(S) = {\partial\Delta/\partial x_4
\over
                (u(x_1,x_2,x_6)u(x_1,x_3,x_5))^{1/2} },$$
where $u$ is the function defined by Equation 8.1.3.

The partial derivative of $\cos\dih(S)$ with respect
to $x_3$ is
$${2x_1 (\partial\Delta/\partial x_2)\over u(x_1,x_2,x_6)^{1/2}
 u(x_1,x_3,x_5)^{3/2}},$$
so that the sign of the partial derivative is determined by the
sign of $\partial\Delta/\partial x_2$.  Similar considerations
apply to the partial derivatives with respect to $x_2$, $x_5$, and
$x_6$.

{\bf Lemma 8.3.2.}
Let $S$ be a quasi-regular tetrahedron.
Then the dihedral angle of $S$ along the any edge is
at most
$\arccos\left({-29003\over 96999}\right) = \dihmax\approx 1.874444$.

{\bf Proof:}  Suppose that $\dih(S)>\pi/2$, so that the numerator $n=
\partial\Delta/\partial x_4$
in Lemma 8.3.1 is negative.  To bound $\cos\dih(S)$ from below, we minimize
$u(x_1,x_2,x_6)$, $u(x_1,x_3,x_5)$, and the numerator $n$ over the variables
$x_2,x_3,x_5$, and $x_6$.

We have $\partial_i u = x_j+x_k-x_i >0$ on the indicated domain,
for $i=2,3,5$, and $6$, so that
$$u(x_1,x_2,x_6) u(x_1,x_3,x_5)\ge u(x_1,2^2,2^2)^2 = (-x_1^2+16 x_1)^2.$$
Similarly $\partial_i n = x_1+x_j-x_k>0$, for $i=2,3,5$, and $6$, so that
$$n(x_1,\ldots,x_6) > n(x_1,2^2,2^2,x_4,2^2,2^2) = x_1(16-x_1-2x_4).$$
Thus,
$$0\ge\cos\dih(S)\ge {16-x_1-2x_4\over 16-x_1} = 1-{2x_4\over 16-x_1}
        \ge 1- {2(2.51)^2\over 16-2.51^2} = {-29003\over 96999}.
        $$\qed

\bigskip
\line{\bf 8.4. The Solid Angles of a Delaunay Simplex\hfill}
\bigskip

Let $S$ be a Delaunay simplex with vertices $v_0,\ldots,v_3$.
By the
solid angle
$\sol_i(S)$ of the simplex at the vertex $i$, we mean 
$3\Vol(S\cap B_i)$,
where ${B}_i$
is a unit ball centered at $v_i$.
For simplicity,
suppose that $v_0$ is located at the origin, and that in the notation of Diagram
8.1.2, the
vertices $v_1,v_2,v_3$ are the edges of lengths $y_1$, $y_2$, and $y_3$.
We write $\sol$ for $\sol_0$.
Set $$a(y_1,y_2,\ldots,y_6) = y_1y_2y_3 +
{1\over2}y_1(y_2^2+y_3^2-y_4^2) + {1\over2}y_2(y_1^2+y_3^2-y_5^2) +
{1\over2}y_3(y_1^2+y_2^2-y_6^2).
\tag 8.4.1$$

{\bf Lemma 8.4.2.}  
$$\sol(S) = {2} \arccot({2 a\over \Delta^{1/2}}).
$$

{\bf Proof:} (See \cite{H2, p.64}).
We use the branch of $\arccot$ taking values in $[0,\pi]$.
\qed

The function $a$ is increasing in $y_1$, $y_2$, and $y_3$ on $[2,4]$
and is decreasing in the variables $y_4$, $y_5$, and $y_6$ on the same interval.

We claim that $a>0$ for any Delaunay simplex.  We will prove this
in the case $\Delta>0$, leaving the degenerate
case $\Delta=0$ to the reader. The area of the face $T$ opposite
the origin is at most $4\sqrt{3}$ (since its edges are all at most
$4$).  There exists a plane that is tangent to the
unit sphere between this face and the unit sphere
\cite{H1,2.1}.  The area
of the radial projection of $T$ to the this plane is less than
the area of a disk $D$ of radius $1.5$ on the plane centered
at the point of tangency.  The solid angle (that is, the area
of the radial projection of $T$ to the unit sphere) is less
than the area of the radial projection of $D$ to the unit sphere.
This area is $2\pi(1-\cos(\arctan(1.5)))<\pi$.  Lemma 8.4.2
now gives the result.  This allows us to use Lemma 8.4.2
in the form $\sol(S) = 2\arctan(\Delta^{1/2}/(2a))$.

{\bf 8.4.3.}  The solid angle is the
area of a spherical triangle.  Let $x$, $y$, and $z$
be the cosines of the radian lengths of the edges of the
triangle. By the spherical law of cosines, the solid angle,
expressed as a function of $x$, $y$, and $z$ is
$$c(x,y,z)+c(y,z,x)+c(z,x,y)-\pi,\quad 
c(x,y,z):=\arccos({x-yz\over \sqrt{(1-y^2)(1-z^2)}}).$$
The partial derivative with respect to $x$ of this expression 
for the solid angle 
is $(-1-x+y+z)/((x+1)\sqrt{t})$, where $t=1-x^2-y^2-z^2+2xyz$.
The second derivative of this expression evaluated at $-1-x+y+z=0$
is $-2(1-y)(1-z)t^{-3/2}\le 0$.  So the 
unique critical point is always a local maximum.
If neither edge at a vertex of a spherical triangle is constrained,
then the triangle can be contracted or expanded by moving the vertex.
  If the lengths of the edges 
are constrained to lie in a product of intervals, then the
minimum area occurs when two of the edges are as short as possible
and the third is at one of the extremes.  The maximum area
is attained when two of the edges are as long as possible and
the third is at one of the extremes or at a critical point:
$x=y+z-1$ (or it symmetries, $y=x+z-1$, $z=x+y-1$).

\bigskip
\line{\bf 8.5. The Compression of a Delaunay Simplex\hfill}
\bigskip

  The compression $\Gamma(S)$ of a Delaunay
simplex $S$ is defined as
$$\Gamma(S) = -\delta_{oct} \Vol(S) + \sum_{i=1}^4\sol_i(S)/3,\tag 8.5.1$$
where $$\delta_{oct} = {-3\pi + 12 \arccos(1/\sqrt{3})\over2\sqrt{2}} \approx
      0.72.\tag 8.5.2$$
We let $a_0=a(y_1,y_2,y_3,y_4,y_5,y_6)$, where $a$ is the function
defined in 8.4.1.  We also let $a_1$, $a_2$, and $a_3$
be the functions $a$ for the vertices denoted $v_1$, $v_2$, and $v_3$ 
in
Diagram 8.1.2.  For example, $a_1=a(y_1,y_5,y_6,y_4,y_2,y_3)$.
Set $t =  \sqrt{\Delta}/2$.
By the results of Sections 8.1 and 8.4, we find that
$\Gamma(S) = -\delta_{oct}\, t/6 +
\sum_0^3 (2/3)\arctan(t/a_i)$.
Recall that $a_i>0$ on a cell $C$.
  We wish to give an elementary upper bound of
$\Gamma$ on $C$.  Set, for $t\ge 0$ and $a=(a_0,a_1,a_2,a_3)\in {\Bbb R_+^4}$,
$$\gamma(t,a) := {-\delta_{oct}\, t\over 6} + \sum_{i=0}^3 {2\over 3}
\arctan({t\over a_i}).$$
The partial derivative of $\gamma$ with respect to $a_i$ is
$-2t/(3(t^2+a_i^2))\le0$, so $\gamma(t,a)\le \gamma(t,a^-)$, where
$a^-=(a_0^-,a_1^-,a_2^-,a_3^-)$
is a lower bound of $a$ on $C$, as   determined
by Section 8.4. 
We study $\gamma(t) = \gamma(t,a^-)$ as a function of $t$ to
obtain an upper bound on $\Gamma(S)$.
We note that
$$\gamma'(t) = {-\delta_{oct}\over 6} + {2\over 3} \sum_{i=0}^3
        {a_i\over (t^2+a_i^2)}\text{ and }
\gamma''(t) = -{4\over 3}\sum {a_i t\over (t^2+a_i^2)^2} \le 0.$$

Upper and lower bounds on $t$ are known from Section 8.1.
An upper bound on $\gamma$ for $t$ in
$[t_{\min},t_{\max}]$ is the maximum of $\ell(t_{\min})$
and $\ell(t_{\max})$, where $\ell$ is the tangent
to $\gamma$ at any point in $[t_{\min},t_{\max}]$.

\bigskip\hbox{}\bigskip
\line{\bf 8.6. Voronoi cells\hfill}
\bigskip

We assume in this section that the simplex 
$S$ has the property that the
circumcenter of each of the three faces with vertex at the
origin lies in the cone at the origin over the face.  This condition is
automatically satisfied if these three faces are acute triangles.
In particular, it is satisfied for a quasi-regular tetrahedron.

When the cone over a
quasi-regular tetrahedron $S$ contains the circumcenter of $S$,
Section 2 sets $\vor(S) = -4\delta_{oct}\Vol(\hat S_0)
	+ 4\sol(S)/3$, where $\hat S_0$ is the intersection of $S$
with a Voronoi cell at the origin.  Otherwise,
$\vor(S)$ is defined as an analytic continuation.
This section gives formulas for $\Vol(\hat S_0)$.

As usual, set $S=S(y_1,\ldots, y_6)$ and $x_i = y_i^2$.
Suppose at first that the circumcenter of $S$ is contained in 
the cone over $S$.
The polyhedron $\hat S_0$ breaks into six pieces, called the
{\it Rogers simplices}.  A
Rogers simplex is the convex hull of the
origin, the midpoint of an edge (the first, second or third edge),
the circumcenter of a face along the 
given edge, and the circumcenter
of $S$.  Each Rogers simplex has the form
$$R=R(a,b,c):= S(a,b,c,(c^2-b^2)^{1/2},(c^2-a^2)^{1/2},(b^2-a^2)^{1/2})$$
for some $1\le a\le b\le c$.  Here $a$ is the half-length of an
edge, $b$ is the circumradius of a face, and $c$ is the
circumradius of the Delaunay simplex.

The volume is $\Vol(R(a,b,c)) = a (b^2-a^2)^{1/2}(c^2-b^2)^{1/2}/6$.
The density $\delta(a,b,c)$ of $R=R(a,b,c)$ is defined as the
ratio of the volume of the intersection of $R$ with a unit ball 
at the origin 
to the volume of $R$.
It follows from the definitions that
$$\vor(S) = \sum 4\Vol(R(a,b,c)) (-\doct+\delta(a,b,c)),
\tag8.6.1$$
where $c$ is the circumradius $\rad(S)$,
and $(a,b)$ runs over the six pairs
$$\align
&(y_1/2,\eta(y_1,y_2,y_6)),\quad (y_2/2,\eta(y_1,y_2,y_6)), \\
&(y_2/2,\eta(y_2,y_3,y_4)),\quad (y_3/2,\eta(y_2,y_3,y_4)), \\
&(y_3/2,\eta(y_3,y_1,y_5)),\quad (y_1/2,\eta(y_3,y_1,y_5)). 
\endalign
$$
Upper and lower bounds on $\Vol(R(a,b,c))$ follow without difficulty
from upper and lower bounds on $a$, $b$, and $c$.
Thus, an upper bound on $\delta(a,b,c)$ leads to an upper bound
on $\vor(S)$.  The next lemma, which is due to Rogers, gives a 
good upper bound \cite{R}.

{\bf Lemma 8.6.2.}  The density $\delta(a,b,c)$ is monotonically
decreasing in each variable for $1<a<b<c$.

{\bf Proof:}  
Let $1\le a_1\le b_1\le c_1$, 
$1\le a_2\le b_2\le c_2$,
$a_2\le a_1$, $b_2\le b_1$, $c_2\le c_1$.  
The points of $R(a_1,b_1,c_1)$ are realized geometrically
by linear combinations
$${\bold s}_1 = 
\lambda_1 (a_1,0,0) + \lambda_2 (a_1,
(b_1^2-a_1^2)^{1/2},0) +
        \lambda_3 (a_1,(b_1^2-a_1^2)^{1/2},(c_1^2-b_1^2)^{1/2}),$$
where $\lambda_1,\lambda_2,\lambda_3\ge 0$ and $\lambda_1+\lambda_2
        +\lambda_3 \le 1$.
The points of $R(a_2,b_2,c_2)$ are realized geometrically
by linear combinations
$${\bold s}_2 = 
\lambda_1 (a_2,0,0) + \lambda_2 (a_2,
(b_2^2-a_2^2)^{1/2},0) +
        \lambda_3 (a_2,(b_2^2-a_2^2)^{1/2},(c_2^2-b_2^2)^{1/2}),$$
with the same restrictions on $\lambda_i$.
Then
$$|{\bold s}_1|^2 - |{\bold s}_2|^2 =
        \lambda_1(\lambda_1+2\lambda_2+2\lambda_3)(a_1^2-a_2^2) +
        \lambda_2(\lambda_2+2\lambda_3)(b_1^2-b_2^2)+
        \lambda_3^2 (c_1^2-c_2^2).$$
So $|{\bold s}_1|^2 \ge |{\bold s}_2|^2$.  This means that
the linear transformation
${\bold s}_1\mapsto {\bold s}_2$
that carries the simplex $S_1$ to $S_2$
moves points of the simplex $S_1$ closer to the origin.
In particular, the linear transformation carries the 
part in $S_1$ of the unit ball at the origin
into the unit ball.
This means that the density of $R(a_1,b_1,c_1)$ is at most
that of $R(a_2,b_2,c_2)$.  \qed

\bigskip
If the circumcenter of $S$ is not in the cone over $S$,
then the analytic continuation gives
$$\vor(S) = \sum_R 4\epsilon_R\Vol(R(a,b,c))(-\doct+\delta(a,b,c)),$$
where $\epsilon_R=1$ if the face of the 
Delaunay simplex $S$ corresponding
to $R$ has positive orientation, and $\epsilon_R=-1$ otherwise.
(The face of $S$ ``corresponding'' to $R(a,b,c)$ is the one
used to compute the circumradius $b$.)

\bigskip
{\bf 8.6.3.}
A calculation based on the explicit coordinates
of $S$ and its circumcenter given in 8.1.4 shows that
$$\epsilon_R \Vol(R(y_1/2,\eta(y_1,y_2,y_6),\rad(S)))=
{ x_1(x_2+x_6-x_1)
	\chi(x_4,x_5,x_3,x_1,x_2,x_6)\over
	48 u(x_1,x_2,x_6)\Delta(x_1,\ldots,x_6)^{1/2}}.$$
If the circumcenter of $S$ is not contained in $S$, then
the same formula holds
by analytic continuation.  
By definition, $\epsilon_R = -1$
exactly when the function $\chi$ is negative.
Although this formula
is more explicit than the earlier formula, it tends
to give weaker estimates of $\vor(S)$ and was not used
in the calculations in Section 9.

\bigskip
{\bf 8.6.4.}
There is another approximation to $\vor(S)$ that
will be useful.
Set $S_y = S(2,2,2,y,y,y)$. (We hope there is no confusion
with the previous notation $S_i$.)
For $1\le a\le b\le c$, let $\Vol(R(a,b,c)) =
 a((b^2-a^2)(c^2-b^2))^{1/2}/6$ be as above.
Set $r(a) = \Vol(R(a,\eta(2,2,2a),1.41))$. 

{\bf Lemma 8.6.5.}  Assume that the circumradius of a quasi-regular
tetrahedron $S$
is at least $1.41$, and that $6\le y_1+y_2+y_3\le 6.3$.  Set
$a= (y_1+y_2+y_3-4)/2$.  Pick $y$ to satisfy $\sol(S_y)=
\sol(S)$.  Then
$$
\vor(S) \le \vor(S_y)-8\doct(1-1/a^3)r(a)
$$

{\bf Proof:}
If $S$ is any quasi-regular tetrahedron, let $S_{\tan}$ be the simplex
defining the ``tangent'' Voronoi cell, that is,
$S_{\tan}$ is the simplex with the same origin
that cuts out the same spherical
triangle as $S$ on the unit sphere, but that satisfies
$y_1=y_2=y_3=2$.  
The lengths of the fourth, fifth, and sixth edges of $S_{\tan}$
are between $\sqrt{8-2(2.3)}=\sqrt{3.4}$ and $2.51$.  The
faces of $S_{\tan}$ are acute triangles.  A calculation
similar to the proof in 8.2.5, 
based on $\chi(3.4,3.4,4,4,4,2.51^2)>0$,
shows that the circumcenter of $S_{\tan}$ is contained in
the cone over $S_{\tan}$.

Since $S_{\tan}$ is obtained by ``truncating'' $S$, we observe that
$\vor(S)=\vor(S_{\tan})-4\doct\Vol(\hat S\backslash 
\hat S_{{\tan}})$,
where $\hat S$ and $\hat S_{{\tan}}$
are the pieces of Voronoi cells denoted $\hat S_0$ in Section
2 for $S$ and $S_{\tan}$, respectively.  
By a convexity result of
L. Fejes T\'oth, $\vor(S_{\tan})\le \vor(S_y)$
\cite{FT,p.125}.
(This inequality relies on the fact that
the circumcenter of $S_{\tan}$ 
is contained in the cone over $S_{\tan}$.)

Let $a_1$ be the half-length of the first, second, or
third edge.
Since $\Vol(R(a_1,b,c))$ is increasing in $c$, we obtain a lower bound
on $\Vol(R(a_1,b,c))$ for $c=1.41$, the
lower bound on the circumradius of $S$.  
The function $\Vol(R(a_1,b,1.41))$,  
considered as a function of $b$, has
at most one critical point in $[a_1,1.41]$ and it is always
a local maximum (by a second derivative test).  
  Thus, 
$$\Vol(R(a_1,b,c))\ge 
   \min(\Vol(R(a_1,b_{\min},1.41)),
   \Vol(R(a_1,b_{\max},1.41))),\tag8.6.6$$
where $b_{\min}$ and $b_{\max}$ are upper and lower bounds
on $b\in[a_1,1.41]$.  A lower bound on $b$ is $\eta(2,2,2a_1)$,
which means that $\Vol(R(a_1,b_{\min},1.41))$ may be replaced
with $r(a_1)$ in the Inequality 8.6.6.
By Heron's formula, the function $\eta(a,b,c)$, for acute
triangles, is convex in pairs of variables:
$$
\eta_{aa}\eta_{bb}-\eta_{ab}^2=
{\eta^6 (a^2+b^2-c^2)(a^2-b^2+c^2)(-a^2+b^2+c^2)(a^2+b^2+c^2)
\over a^6b^6c^4}>0.$$
Thus, 
an upper bound on $b_{\max}$ is $\eta(2,2.51,2a)$, where
$a = (y_1+y_2+y_3-4)/2$.  This means
that $$\Vol(R((a_1,b_{\max},1.41)))$$ may be replaced with
the function
$\Vol(R(a_1,\eta(2,2.51,2a),1.41))$ in the Inequality 8.6.6.
Now $1\le a_1\le a$, and $\Vol(R(a_1,b,c))$ 
is decreasing in the first variable
(for $1\le a_1$ and $b\le \sqrt{2}$), so we may use
the lower bound $$\Vol(R(a,\eta(2,2.51,2a),1.41))$$ instead.
By  Calculations 9.20.2 and 9.20.3,
we conclude that $\Vol(R(a_1,b,c))\ge r(a)$.
This lower bound is valid
for each Rogers simplex.  The volume of 
$\hat S\backslash \hat S_{{\tan}}$
is then at least
$$((1-8/y_1^3) + (1-8/y_2^3) + (1-8/y_3^3))2r(a).
$$
The concavity
of $1-8/y^3$ gives $\sum_{i=1}^3 
(1- 8/y_i^3)\ge 1-1/a^3$.
We have established that 
$$\Vol(\hat S\backslash \hat S_{{\tan}})\ge (1-1/a^3) 2
r(a).$$
The result follows.  \qed

{\bf 8.6.7.}  We conclude our discussion of Voronoi cells with
a few additional comments about the case in which analytic
continuation is used to define $\vor(S)$, with $S$  a
quasi-regular tetrahedron.  Assume the circumcenter $c$ of $S$
lies outside $S$ and that the face $T$ of $S$ with negative
orientation is the one bounded by the first, second, and sixth
edges.  It follows from Section 3 that $y_1, y_2, y_6\in[2.3,2.51]$
and $y_3, y_4, y_5\in[2,2.15]$.

Let $p_1$ (resp. $p_2$) be the point on $T$ equidistant from
the origin, $v_3$, and $v_1$ (resp. $v_2$).
$$|p_1|^2 = 
{x_1\over 4} + {x_1 u(x_1,x_2,x_6)(-x_1+x_3+x_5)^2\over
4 (\partial\Delta/\partial x_4)^2}.$$
Let $p_0$ be the circumcenter of $T$ (see Diagram 8.6.8).

\gram|2|8.6.8|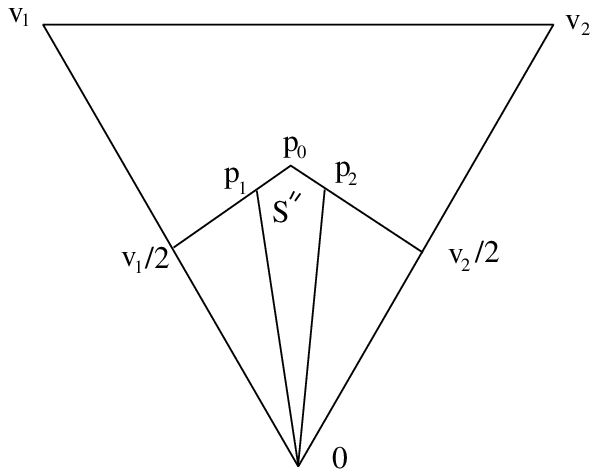|  
$\epsilon_R=-1$ for the Rogers simplex with vertices the origin,
$p_0$, $c$ and $v_1/2$.  It lies outside $S$.  The other Rogers
simplex along the first edge has $\epsilon_R=1$, so that the
part common to both of these Rogers simplices cancels in the
definition of $\vor(S)$. This means that $\vor(S)$ becomes the
sum of the usual contributions from the two Rogers simplices
along the third edge,
$$4\Vol(R(a,b,c))(-\doct+\delta(a,b,c))$$
for $(a,b,c) = (y_1/2,\eta(y_1,y_3,y_5),|p_1|))$ and
  $(y_2/2,\eta(y_2,y_3,y_4),|p_2|)$,
and
$$v(S''):= 4\doct\Vol(S'') - 4\sol(S'')/3,$$
where $S''$ is the convex hull of $0$, $p_0$, $p_1$, $p_2$, and $c$.

We claim that in any cell satisfying the constraints given above,
$|p_1|$ is minimized by making
$y_1$, $y_2$, $y_6$ as large as possible and $y_3$, $y_4$, $y_5$
as small as possible.  To show this, one has to write out the
derivatives explicitly from the formula for $|p_1|^2$ given
above.  The partial derivatives 
with respect to  $x_2,\ldots,x_6$ factor into products
of the polynomials 
$\partial\Delta/\partial x_i$, 
for $i=3,4,5$, $x_1$, $(x_3+x_5-x_1)$,
$u(x_1,x_2,x_6)$, $(-x_1x_4+x_2x_5+x_1x_6-x_5x_6)$, 
and $(x_1x_2-x_2x_3-x_1x_4+x_3x_6)$.
The signs of these polynomials are easily determined.  The
partial with respect to $x_1$ is complicated (the numerator
has 88 terms), and we had to
resort to the method of subdivision to determine its sign.
We omit the details.

To complete our estimate, we describe an upper bound on $v(S'')$.
The ratio $\sol(S'')/(3\Vol(S''))$ is at least $1/\rad(S)^3$, because
$S''$ is contained in a sphere of radius $\rad(S)$, centered at
the origin. We have
$${\Vol(S'')\over2}\le {1\over 6} |p_0-p_1|\,|c-p_0|\,|p_0|,$$
because the convex hull of $0,p_0,p_1$, and $c$, which is one side
of $S''$, is a pyramid with base the right triangle $(p_1,p_0,c)$
and height at most $|p_0|=\eta(y_1,y_2,y_6)$.  Of course,
$|c-p_0|^2 = \rad(S)^2-\eta(y_1,y_2,y_6)^2$.  We now have a bound
on $v(S'')$ in terms of quantities that have been studied in
Section 8, if we rely on the bound $|p_0-p_1|\le 0.1381$.
Write $p_1=p_1(y_1,\ldots,y_6)$.  
This bound is obtained from the following inequalities.
$$\align
 |p_0-p_1| &= (\eta(y_1,y_2,y_6)^2-x_1/4)^{1/2} - 
     (|p_1|^2-x_1/4)^{1/2}\\
  &\le (\eta(y_1,2.51,2.51)^2-x_1/4)^{1/2} - 
      (|p_1(y_1,2.51,2,2,2,2.51)|^2-x_1/4)^{1/2}\\
  &\le (\eta(2.51,2.51,2.51)^2-2.51^2/4)^{1/2} - \\
      &\qquad\qquad(|p_1(2.51,2.51,2,2,2,2.51)|^2-2.51^2/4)^{1/2}\\
  & < 0.1381.
\endalign
$$
The inequality that replaces $y_1$ with $2.51$ results from
Calculation 9.21.

\bigskip
\line{\bf 8.7. A Final Reduction\hfill}
\bigskip

Let $S$ be a Delaunay simplex.
Suppose
that the lengths of the edges $y_1$, $y_5$, and $y_6$ are greater than
$2$.  Let $S'$ be a simplex formed by contracting the vertex joining
edges $1$, $5$, and $6$ along the first edge by a small amount.  We assume
that the lengths  $y'_1$, $y'_5$, and $y'_6$ of the new edges
are still at
least $2$ and that the circumradius of $S'$ is at most 2, so
that $S'$ is a Delaunay simplex.

\bigskip
{\bf Proposition 8.7.1.}  $\Gamma(S') > \Gamma(S)$.

\bigskip
We write $\sol_i$, for $i=1,2,3$, for the solid angles at the three
vertices $p_1$, $p_2$, and $p_3$ of $S$ terminating
the edges $1$, $2$, and $3$. Let $p_1'$ be the vertex terminating
edge $1$ of $S'$.
Similarly,
we write $\sol'_i$, for $i=1,2,3$, for the solid angles at the
corresponding vertices of $S'$.  We set $\Vol(V) = \Vol(S)-\Vol(S')$ and
$w_i = \sol_i-\sol'_i$.  It follows directly from the
construction of $S'$ that $w_2$ and $w_3$ are positive.
The dihedral angle $\alpha$ along the first edge is the same for
$S$ and $S'$. The angle $\beta_i$ of the triangle $(0,p_1,p_i)$ 
at $p_1$ is less than the angle $\beta'_i$ of the triangle
$(0,p'_1,p_i)$ at $p'_1$, for $i=2,3$.  It follows that $w_1$
is negative, since $-w_1$ is the area of the quadrilateral
region of Diagram 8.7.2 on the unit sphere.

\gram|2|8.7.2|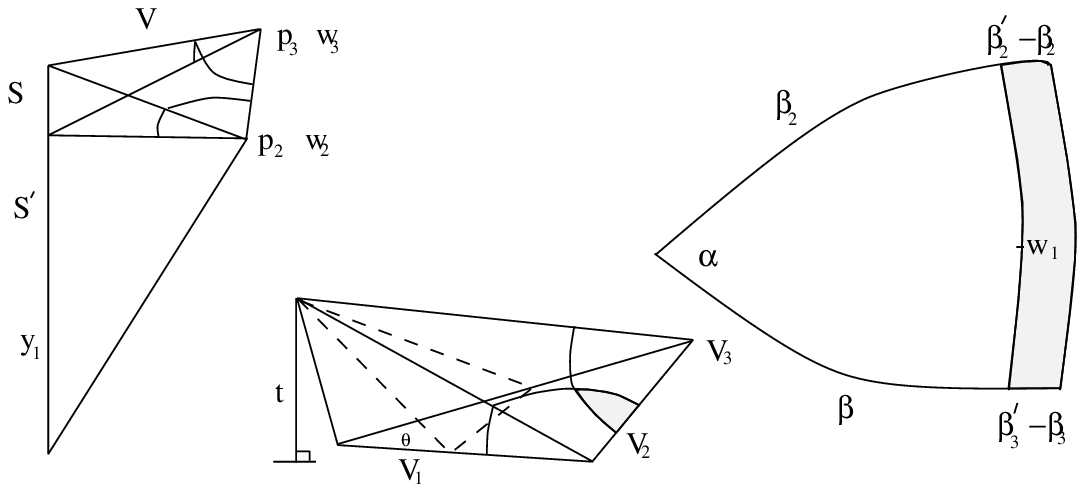|

\bigskip
{\bf Lemma 8.7.3.}  $\delta_{oct} \Vol(V) > w_2/3 + w_3/3 $.

\bigskip
The lemma immediately implies the proposition because $w_1<0$
and
$$\Gamma(S') -\Gamma(S) = -w_1/3 -w_2/3 -w_3/3 + \delta_{oct}\Vol(V).$$

\bigskip
{\bf Proof:}  
Let $T'$ be the face of $S'$ with
vertices $p'_1$, $p_2$, and $p_3$.
We consider $S'$ as a function of $t$, where $t$ is the
distance from $p_1$ to the plane containing $T'$.
(See Diagram 8.7.2.)  
It is enough to establish the lemma
for $t$ infinitesimal.

As shown in Diagram 8.7.2,  let $V_1$ be the pyramid formed
by intersecting $V$ with the plane through $p_1$ that meets
the fifth edge at distance $t_0=1.15$ from $p_3$ and the sixth
edge at distance $t_0$ from $p_2$.
Also, let the intersection of $V$ with
a ball of radius $t_0$ centered at $p_i$ be denoted $V_i$, for
$i=2,3$.  For $t$ sufficiently small, the region $V_1$ is (essentially)
disjoint from $V_2$ and $V_3$.

We claim that $\Vol(V_1) > \Vol(V_2\cap V_3)$.  
 Let $\theta$ be the angle of $T'$
subtended by the fifth and sixth edges of $S'$.  Then
$\Vol(V_1) = B t/3$, where $B$ is the area of the intersection
of $V_1$ and $T'$: 
$$B = \sin\theta (y_5'-t_0)(y_6'-t_0)/2.\tag 8.7.4$$
Since $S'$ is a Delaunay simplex, the estimates $\pi/6\le\theta\le2\pi/3$
from \cite{H1,2.3} hold.  In particular, $\sin\theta\ge 0.5$,
so $B\ge 0.25(2-1.15)^2$, and $\Vol(V_1)> 0.06\,t$.

If $\Vol(V_2\cap V_3)$ is nonempty, the fourth edge of $S$ must
have length less than $2t_0$.  The dihedral angle $\alpha'$ of $V$ along
the fourth edge is then less than the tangent of the angle, which
is at most $t+ {O}(t^2)$, in Landau's notation.  
As in \cite{H1,5}, we obtain the
estimate
$$\Vol(V_2\cap V_3) \le {\alpha'} \int_1^{t_0}
        t_0^2-t^2\,dt = { (t_0-1)^2(2t_0+1)\alpha'\over 3}
        < 0.025\, t + {O}(t^2).$$
This establishes the claim.

Thus, for $t$ sufficiently small $$\align
  \delta_{oct}\Vol(V) &\ge \delta_{oct}(\Vol(V_1) +\Vol(V_2)
        +\Vol(V_3) - \Vol(V_2\cap V_3)) \\
        &>
        (\delta_{oct} t_0^3) \left( {w_2\over3} + {w_3\over3}\right)
        > 1.09 \left( {w_2\over3} + {w_3\over3}\right).
\endalign$$
\qed

\bigskip
\centerline{\bf Section 9. Floating-Point Calculations}
\bigskip

This section describes various inequalities that have been established
by the method of subdivision on SUN workstations.  The full source code
(in {\tt C++}) for these calculations is available \cite{H6}.

Floating-point operations on computers are subject to round-off errors,
making many machine computations unreliable.  Methods of
interval arithmetic give users control
over round-off errors \cite{Int}.  These methods may be reliably
implemented on machines that allow arithmetic
with directed rounding,
for example those conforming to the IEEE/ANSI standard 754 \cite{W},\cite{IEEE},
\cite{NR}.

Interval arithmetic produces an interval in the real line
that is guaranteed to contain the result of an arithmetic
operation.  As the round-off errors accumulate, the interval
grows wider, and the correct answer remains trapped in the
interval.
Apart from the risk of compiler errors and defective hardware,
a bound established by interval arithmetic is as reliable
as a result established by integer arithmetic on a computer.
We have used interval arithmetic wherever computer precision
is a potential issue (Calculations 9.1 -- 9.19, in particular).

Every inequality of this section has been reduced to a finite
number of inequalities of the form $r(x_0)<0$, where $r$ is
a rational function of $x\in {\Bbb R}^n$ and $x_0$ is a given
element in ${\Bbb R}^n$.  To evaluate each rational expression,
interval arithmetic is used to obtain an interval $Y$ containing
$r(x_0)$.  The stronger inequality,
$y<0$ for all $y\in Y$, which may be verified by computer,
implies that $r(x_0)<0$.

To reduce the calculations to rational expressions $r(x_0)$,
rational approximations to the functions $\sqrt{x}$, $\arctan(x)$,
and $\arccos(x)$ with explicit error bounds are required.
These were obtained from \cite{CA}.  Reliable approximations to
various constants (such as $\pi$, $\sqrt{2}$, and $\doct$) with
explicit error bounds are also required.  These were obtained
in {\it Mathematica\/} and were double checked against {\it Maple}.

\bigskip
Let $S=S(y_1,\ldots,y_6)$ be a quasi-regular 
tetrahedron.
We label the indices as in
Diagram 8.1.2.  Let $\Gamma=\Gamma(S(y_1,\ldots,y_6))$ be the compression.
We also
let $\dih=\dih(S(y_1,\ldots,y_6))$ be the dihedral angle along the
first edge,
and let $\sol=\sol(S(y_1,\ldots,y_6))$ be the solid angle at the
origin, that is, the solid angle formed by the first, second, and third
edges ($y_1,y_2,y_3$) of $S$.

All of the following inequalities are to be considered as inequalities
of analytic functions of $y_1,\ldots,y_6$.  Although each of the
calculations is expressed as an inequality between functions of
six variables, Lemma 8.7.1 has been invoked repeatedly to reduce
the number of variables to three or four.
For instance, suppose that we wish to establish $I(\sol,\Gamma)<0$,
where $I(\sol,\Gamma)$ is an expression in $\Gamma(S(y_1,\ldots,y_6))$
and $\sol(S(y_1,\ldots,y_6))$.  Invoking 8.7.1 three times, we may
assume that the vertices marked $v_1$, $v_2$, and $v_3$ 
in Diagram 8.1.2
each terminate an edge of minimal length.  It is then sufficient
to establish the inequality in seven situations of
smaller dimension; that is, we may assume the  edges
$i\in I$ have minimal length, where $I$
is one of
$$\{1,4\},\ \{2,5\},\ \{3,6\},\ \{4,5\},\ \{4,6\},\ \{5,6\},\ \{1,2,3\}.$$
Similarly, for an inequality in $\Gamma$ and dihedral angle, we reduce to
the seven cases
$$I=\{1,4\},\ \{2,5\},\ \{3,6\},\ \{1,2,3\},\ \{1,5,6\},\
        \{3,4,5\},\ \{2,4,6\}.$$

The first two calculations
are inequalities of the compression $\Gamma$ of a quasi-regular tetrahedron.

\smallskip
{\bf Calculation 9.1.}  $\Gamma\le 1\,\pt$. \quad 

\medskip
This first inequality and Calculation 9.3 are the only ones that
are not strict inequalities.  Set $S_0 = S(2,2,2,2,2,2)$.
By definition, $\Gamma(S_0) = 1\,\pt$.
We must give a direct proof that $S_0$ gives the maximum
in an explicit neighborhood.  Then we use the method of subdivision
to bound $\Gamma$ away from $1\,\pt$ outside the given neighborhood.
An infinitesimal version of the following result is
proved in \cite{H1}.

\bigskip
{\bf Lemma 9.1.1.}   If $y_i\in [2,2.06]$, for $i=1,\ldots,6$, then
$\Gamma(S(y_1,\ldots,y_6)) \le \Gamma(S_0)$, with equality if and
only if $S(y_1,\ldots,y_6) = S_0$.

\bigskip
Set $a_{00}=a(2,2,2,2,2,2) = 20$, $\Delta_0=\Delta(2^2,2^2,2^2,2^2,2^2,2^2)$,
$t_0=\sqrt{\Delta_0}/2 = 4\sqrt{2}$, and
$b_0=(2/3)(1+t_0^2/a_{00}^2)^{-1} = 50/81$.
Set $f=\max_i(y_i-2)\le 0.06$
and $a^-=a(2,2,2,2+f,2+f,2+f)
= 20-12f-3f^2$.
As in Section 8, we set $t=\sqrt{\Delta}/2$.
We have $\arctan(x)\le \arctan(x_0) + (x-x_0)/(1+x_0^2)$,
if $x,x_0\ge0$. We will verify below that $\Delta\ge \Delta_0$
under the restrictions given above.  Then
$$
\align
\Gamma(S) &= -{\delta_{oct} t\over 6} + {2\over3} \sum_{i=0}^3\arctan(t/a_i)\\
        &\le \Gamma(S_{0}) -{\delta_{oct}(t-t_0)\over 6} +
                b_0 \sum_{i=0}^3
                \left({t\over a_i} - {t_0\over a_{00}}\right)\\
        &=\Gamma(S_{0}) - {\delta_{oct}(\Delta-\Delta_0)\over 24(t+t_0)}
                +b_0\sum_{i=0}^3 {t-t_0\over a_i}\\
        &\qquad\quad
                +{b_0 t_0\over a_{00}^2}\sum_{i=0}^3 {(a_{00}-a_i)}
                +{b_0t_0\over a^2_0}\sum_{i=0}^3 {(a_{00}-a_i)^2\over a_i}\\
        &\le \Gamma(S_{0}) + {(\Delta-\Delta_0)\over t+t_0}\left(
                -{\delta_{oct}\over 24 } + {b_0\over a^-}
                \right) +
                {b_0t_0\over a_{00}^2}\sum_{i=0}^3 (a_{00}-a_i)
                \\
        &\qquad\quad +
                {b_0t_0\over a^2_0a^-}\sum_{i=0}^3 (a_{00}-a_i)^2.\tag9.1.2
\endalign
$$
Set $c_0 = -\delta_{oct}/24 + b_0/a^- > 0$.
Write $y_i = 2+f_i$, with $0\le f_i\le f$.  Set $x_i = 4+e_i =
        y_i^2 = 4+4f_i+f_i^2$.  Then $f_i\le e_i/4$, for $i=1,\ldots,6$.
Set $e=4f+f^2$.  Set $\tilde\Delta(e_1,\ldots,e_6)=\Delta(x_1,\ldots,x_6)$.
We find that $\partial^2\tilde\Delta/\partial e_2\partial e_1
        =4 +e_4+e_5-e_6>0$ and similarly that
$\partial^2\tilde\Delta/\partial e_i\partial e_1>0$, for $i=3,5,6$.
So $\partial\tilde\Delta/\partial e_1$
 is minimized by taking $e_2=e_3=e_5=e_6=0$,
and this partial derivative is at least
$${\partial\tilde\Delta\over\partial e_1}
        (e_1,0,0,e_4,0,0) = 16-8e_1-2e_1e_4-e_4^2
        \ge 16-8e-3e^2 >0.$$
Thus $\Delta\ge\Delta_0$.
The partial derivative $\partial\tilde\Delta/\partial e_1$ is at most
$$\align
        {\partial\tilde\Delta\over \partial e_1} (e_1,e,e,e_4,e,e) &=
        16+16e - 8e_1 + 4 e e_4 - 2e_1 e_4-e_4^2 \le 16+16 e + 4 e e_4 - e_4^2 \\        &\le 16+16e + 3 e^2.
\endalign
$$
So $\Delta-\Delta_0\le (16+16e+3e^2)(e_1+\cdots+e_6)$.

We expand $\sum_{i=0}^3 (a_i-a_{00}) = 5(e_1+\cdots+e_6)+h_1+h_2+h_3$ as a sum of
homogeneous polynomials $h_i$ of degree $i$ in $f_1,\ldots,f_6$.  A
calculation shows that $h_1=0$.  Also $h_2$ is quadratic in each variable $f_i$
with negative leading coefficient, so $h_2$ attains its minimum at an extreme
point of the cube $[0,f]^6$.  A calculation then shows that $h_2\ge -14f^2$ on
$[0,f]^6$.  By discarding all the positive terms of $h_3$, we find that
$$h_3 \ge -f_1^2 f_4 - f_1 f_4^2 - f_2^2 f_5 -f_2 f_5^2 - f_3^2 f_6 - f_3 f_6^2 \ge -6f^3.$$
Thus, $\sum_{i=0}^3(a_i-a_{00}) \ge 5(e_1+\cdots+e_6) - 14 f^2 -6f^3$.  Since $e = \max(e_i)$,
we find that $f\le e/4 \le (e_1+\cdots+e_6)/4$.  This gives
$$
\sum_{i=0}^3 (a_{00}-a_i) \le (e_1+\cdots+e_6)(-5 + {7f\over2} + {3f^2\over 2}).$$

Similarly, we expand $\sum_{i=0}^3(a_{00}-a_i)^2$ as a polynomial in $f_1,\ldots,f_6$.
To obtain an upper bound, we discard all the negative terms of the polynomial and
evaluate all the positive terms at $(f_1,\ldots,f_6) = (f,\ldots,f)$.  This gives$$\align
\sum_{i=0}^3 (a_{00}-a_i)^2 &\le
        4944 f^2 + 7296 f^3 + 3684 f^4 + 828 f^5 + 73 f^6\\
        &\le {1\over 4}(e_1+\cdots+e_6)(4944 f + 7296 f^2 + 3684 f^3 + 828 f^4 +
73 f^5).
\endalign
$$
We now insert these estimates back into Inequality 9.1.2.  This gives
$$
\align
{ \Gamma(S)-\Gamma(S_0)\over (e_1+\cdots+e_6)} &\le
        {c_0\over 2t_0}(16+16e+3e^2) + {b_0t_0\over a_{00}^2} (-5+3.5 f + 1.5f^2)\\
        &\qquad\quad +
                {b_0 t_0\over 4 a_{00}^2 a^-}
                (4944 f + 7296 f^2 + 3684 f^3 + 828 f^4 + 73 f^5).
\tag9.1.3
\endalign
$$
The right-hand side of this inequality is a rational function of $f$.  (Both
$a^-$ and $e$ depend on $f$.)  Each of the three terms on the right-hand side is
increasing in $f$.  Therefore, the right-hand side reaches its maximum 
at $f=0.06$.
Direct evaluation at $f=0.06$ gives $\Gamma(S)-\Gamma(S_0) \le -0.00156(e_1+\cdots+e_6)$.
\qed

To verify Calculation 9.14, the computer examined only $7$ cells.
Calculation 9.1 required
1,899 cells.  
Calculation 9.6.1
required over 2 million cells.  The number of cells required in
the verification of the other inequalities falls between these extremes.

\bigskip
{\bf Calculation 9.2.} $\Gamma < 0.5\,\pt$, if $y_1\in [2.2,2.51]$.

\bigskip
The next several calculations are concerned with the relationship
between the dihedral angle and the compression.

\smallskip
{\bf Calculation 9.3.}  $$\dih(S(y_1,\ldots,y_6))\ge 
\dihmin := \dih(S(2,2.51,2,2,2.51,2))
\approx 0.8639.$$

Since this bound is realized by a simplex, we must carry out the appropriate
local analysis in a neighborhood of $S(2,2.51,2,2,2.51,2)$.

\bigskip
{\bf Lemma 9.3.1.}  Suppose that $2\le y_i\le 2.2$, for $i=1,3,6$,
and that $2\le y_i\le 2.51$, for $i=2,4,5$.  Then
$\dih(S(y_1,\ldots,y_6))\ge \dihmin$.

\bigskip
{\bf Proof:}  This is an application of Lemma 8.3.1.  By that
lemma, $\dih(S)$ is increasing in $x_4$, so we fix $x_4=4$.  
The sign of $\partial\cos\dih/\partial x_2$ is the sign of
$\partial\Delta/\partial x_3$, and a simple estimate based on the
explicit formulas of Section 8 shows that $\partial\Delta/\partial x_3>0$
under the given constraints.
Thus, we minimize $\dih(S)$ by setting $x_2=2.51^2$.  By symmetry,
we set $x_5=2.51^2$.

Now consider $\dih(S)$ as a function of $x_1,x_3$, and $x_6$.  
The sign of $\partial\cos\dih/\partial x_3$ is the sign of
$\partial\Delta(x_1,2.51^2,x_3,2^2,2.51^2,x_6)/\partial x_2$.
The maximum of this partial 
(about $-2.39593$) is attained when $x_1$, $x_3$, and $x_6$
are as large as possible: $x_1=x_3=x_6=2.2^2$.
So $\dih(S)$ is increasing in $x_3$.  We take $x_3=4$, and by symmetry
$x_6=4$.

We have
$$
{\partial \cos\dih(S)\over \partial x_1} = {2t_1(x_1)\over
                u(x_1,x_2,x_6)^{3/2} u(x_1,x_3,x_5)^{3/2}},$$
where $t_1(x_1)\approx 464.622-1865.14x_1+326.954x_1^2-92.9817x_1^3+4x_1^4$.
An estimate of the derivative of $t_1(x)$ shows that $t_1(x)$
attains its maximum at $x_1=4$,
and $t_1(4)< -6691<0$.  Thus $\dih(S)$ is minimized when
$x_1=4$.\qed

\bigskip

\smallskip
{\bf Calculation 9.4.}  $\Gamma < 0.378979\dih -0.410894$.

\smallskip
{\bf Lemma 9.5.}  If $S_1$, $S_2$, $S_3$, and $S_4$ are any four
tetrahedra such that $\dih(S_1) + \dih(S_2) + \dih(S_3) +\dih(S_4) \ge 2\pi$,
then $\Gamma(S_1)+\cdots+\Gamma(S_4) < 0.33\,\pt$.

\smallskip
This lemma is a consequence
of the following three calculations, each established by interval arithmetic
and the method of subdivision.

\smallskip
{\bf Calculation 9.5.1.}  $$\Gamma < -0.19145\dih+ 0.2910494,$$
        provided $\dih\in [\dih(S(2,2,2,2,2,2)),1.42068]$.

\smallskip
{\bf Calculation 9.5.2.}  $$\Gamma < -0.0965385\dih+ 0.1562106,$$
        provided $\dih\in [1.42068,\dih(S(2,2,2,2.51,2,2))]$.

\smallskip
{\bf Calculation 9.5.3.}  $$\Gamma < -0.19145\dih+ 0.31004,$$
        provided $\dih\ge \dih(S(2,2,2,2.51,2,2))$.

To deduce Lemma 9.5, we consider the piecewise linear bound $\ell$ on
$\Gamma$ obtained from these estimates.  The linear
pieces are $\ell_1(x)\le 1\,\pt$,
$\ell_2$, $\ell_3$, and
$\ell_4$ on $[\dihmin,d_1]$, $[d_1,d_2]$, $[d_2,d_3]$, 
and $[d_3,\dihmax]$ , where $d_1= \dih(S(2,2,2,2,2,2))$,
$d_2=1.42068$, and
$d_3 = \dih(S(2,2,2,2.51,2,2))$.
(See Calculation 9.1, and Lemma 8.3.2.)
Diagram 9.5.4 illustrates these linear bounds.  (There are small
discontinuities at $d_1$, $d_2$, and $d_3$ that may be eliminated by
replacing $\ell_i$ by $\ell_i+\epsilon_i$, for some $\epsilon_i>0$, for
$i=1,2$, and $4$.)
 We then ask for the maximum of
$\ell(t_1) + \ell(t_2) + \ell(t_3) + \ell(t_4)$ under the constraint
$(t_1+t_2+t_3+t_4)/4\ge \pi/2 \in [d_2,d_3]$.
Since $\ell(x)$ is constant on $(\dihmin,d_1)$ and decreasing on
$[d_1,\dihmax]$, we may assume
that $t_i\ge d_1$, for all  $i$.  Since the slope of $\ell_2$ is
equal to the slope of $\ell_4$, we may assume that $t_i \ge d_2$,
for all $i$, or that $t_i\le d_3$, for all $i$.
If $t_i\le d_3$, for all $i$, then we find that $t_i\ge 2\pi-3 d_3 > d_2$.
So in either case, $t_i\ge d_2$, for all $i$.
By convexity, an upper bound is
$4\ell_3(\pi/2)< 0.33\,\pt$.

\gram|2|9.5.4|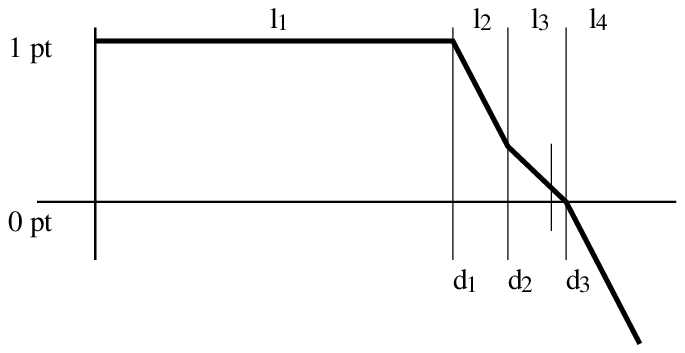|  

\bigskip
{\bf Lemma 9.6.}  If $S_1,\ldots, S_5$ are any five tetrahedra
such that $\dih(S_1)+\cdots+\dih(S_5) \ge 2\pi$, then
$$\Gamma(S_1)+\cdots+\Gamma(S_5) < 4.52\,\pt.$$

\smallskip
This is a consequence of two other calculations.

\smallskip
{\bf Calculation 9.6.1.}  $\Gamma< -0.207045\dih + 0.31023815$,
        provided $\dih\in [d_0,2\pi-4d_0]$, where
        $d_0 =\dih(S(2,2,2,2,2,2))$.

\smallskip
{\bf Calculation 9.6.2.}  $\Gamma < 0.028792018$, if $\dih> 2\pi-4d_0$.

\smallskip
Calculations 9.1, 9.6.1, and 9.6.2 give a piecewise linear bound $\ell$
on $\Gamma$ as a function of dihedral angle.  See Diagram 9.6.3.
(Again, there are minute discontinuities that may be eliminated in the
same manner as before.)
  We claim that $4.52\,\pt >
\ell(t_1)+\cdots+\ell(t_5)$ whenever $t_1+\cdots+t_5\ge 2\pi$.  As in
Lemma 9.5, we may assume that $t_i\ge d_0$.  Since $\ell$ 
is decreasing on $[d_0,\dihmax]$, we may assume that $t_1+\cdots+t_5=2\pi$.
Thus, only the interval $[d_0,2\pi-4d_0]$ is relevant for
the optimization.  On this interval, the bound is linear, so
$\ell(t_1)+\cdots+\ell(t_5) \le 5 \ell(2\pi/5)   < 4.52\,\pt$.

\gram|2|9.6.3|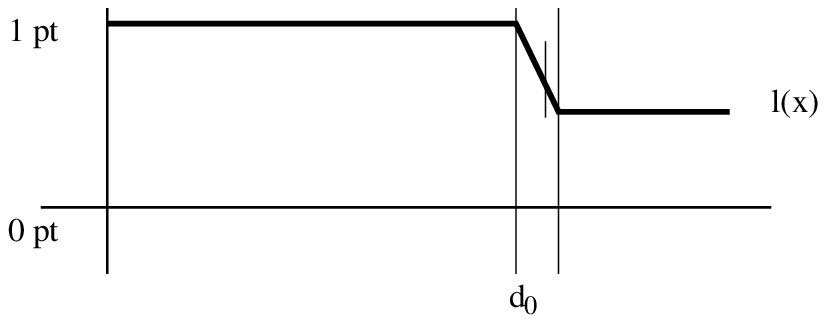|  

\bigskip
{\bf Calculation 9.7.}  $\Gamma < 0.389195\dih -0.435643$,
        if $y_1\ge 2.05$.

\bigskip
The next inequalities relate the solid angles to the
compression. 

\smallskip
{\bf Calculation 9.8.} $\Gamma< -0.37642101\sol + 0.287389$.

\smallskip
{\bf Calculation 9.9.} $\Gamma< 0.446634 \sol -0.190249$.

\smallskip
{\bf Calculation 9.10.}  $\Gamma < -0.419351\sol + 0.2856354 + 0.001$,
if $y_i\in [2,2.1]$, for $i=1,\ldots,6$.

\smallskip
The following calculation involves the circumradius.  We leave
it to the reader to check that the dimension-reduction techniques
of Lemma 8.7.1 may still be applied.

\smallskip
{\bf Calculation 9.11.}  $\Gamma < -0.419351\sol + 0.2856354$,
provided that $y_4, y_5, y_6 \ge 2.1$ and that the circumradius of $S$
is at most $1.41$.

\smallskip
{\bf Calculation 9.12.}  $\Gamma < -0.419351\sol + 0.2856354 - 5(0.001)$,
if $y_i> 2.1$ for some $i$, 
and $y_4\in [2,2.1]$.

\smallskip
{\bf Calculation 9.13.}  $\Gamma < -0.65557\sol + 0.418$,
	if $y_i\in [2,2.1]$, for $i=1,\ldots,6$.

\smallskip
{\bf Calculation 9.14.}  $\sol(S) > 0.21$.

\smallskip
{\bf Calculation 9.15.}  $\Gamma+(K-\sol)/3 < 0.564978\dih  -0.614725$,
where $K = (4\pi-6.48)/12$,
        provided $y_1\in [2,2.05]$, and
        $y_2,y_3\in [2,2.2]$.

\smallskip
{\bf Calculation 9.16.}  $\dih> 0.98$ and $\sol> 0.45$,
provided $y_1\in[2,2.05]$, and $y_2,y_3\in [2,2.2]$.

\bigskip
Let $\Gamma(S)$ be replaced by $\vor(S)$ in each of the Calculations
$9.1$ -- $9.16$ to obtain a new list of inequalities $9.1'$-- $9.16'$.  
(In
$9.11'$, we drop the constraint on the circumradius of $S$.)
We claim that all of the inequalities $9.*'$ hold whenever
$S$ is a quasi-regular tetrahedron of circumradius at least $1.41$.
In fact,  Inequalities $9.1'$, $9.2'$, $9.4'$, $9.5.1'$, $9.5.2'$,
$9.5.3'$, $9.6.1'$, $9.6.2'$, and $9.7'$ follow directly from
$9.17$ and the inequalities $\dihmin\le\dih\le \dihmax$.
Calculations $9.3$, $9.14$, and $9.16$ are independent of $\Gamma$,
and so do not require modification.  Inequalities $9.8'$, $9.11'$, and
$9.12'$ also rely on $9.17$, $9.18$ and $9.19$,  inequality $9.9'$ on
$9.14$, and inequality $9.15'$ on $9.16$.  Inequalities
$9.10'$ and $9.13'$ are vacuous by the comments of 8.2.4.

Write $S=S(y_1,y_2,y_3,y_4,y_5,y_6)$,
$\vor = \vor(S)$, and let $\rad=\rad(S)$ be the circumradius of $S$.

{\bf Lemma 9.17.}  If the circumradius is at least $1.41$,
then $\vor< -1.8\,\pt$.

{\bf Proof:}  If $\sol\ge 0.91882$, then the lemma is a consequence of
Lemma 9.18.  (The proof of Lemma 9.18, under the restriction
$\sol\ge 0.91882$, is independent of the proof of this lemma.)
Assume that $\sol< 0.91882$.

If $S'$ and $S$ are Delaunay stars, related as in Section
8.7, then $\vor(S')> \vor(S)$ because  $\hat S'_0$ is
obtained by slicing a slab from  $\hat S_0$.  This
means that we may apply the dimension-reduction techniques
described at the beginning of this section, unless the
deformation decreases the circumradius to $1.41$.
An interval calculation
establishes the result when $y_1+y_2+y_3\ge 6.3$
and the circumradius constraint is met
(Calculation 9.17.1). 
When
the dimension reduction techniques apply, the dimension
of the search space may be reduced to four, and an interval
calculation similar to the others in this section gives the
result (Calculation 9.17.2).  

We lacked the  computer resources to perform the
interval analysis directly when the dimension-reduction techniques
fail and found it necessary to break the problem up into
smaller pieces when $y_1+y_2+y_3\le6.3$.
We will make use of the following calculations.

{\bf Calculation 9.17.1.}  $\vor(S)<-1.8\,\pt$ provided 
$y_1+y_2+y_3\ge 6.3$, $\rad(S)=1.41$, and $\sol(S)\le 0.91882$.

{\bf Calculation 9.17.2.}  $\vor(S)<-1.8\,\pt$\/ provided that
$\rad(S)>1.41$, 
$y_1+y_2+y_3\ge 6.3$,
$\sol(S)\le 0.912882$, and $S$ lies in one of
the seven subspaces of smaller dimension associated with $I(\sol,\Gamma)$
at the beginning of the section.

{\bf Calculation 9.17.3.1.}  Assume that $y_1+y_2+y_3\le 6.3$
and that $\rad\ge1.41$. Then $\sol\ge 0.767$.

{\bf Calculation 9.17.3.2.} Assume that $y_1+y_2+y_3\le 6.192$
and that $\rad\ge1.41$.  Then $\sol\ge 0.83$.

{\bf Calculation 9.17.3.3.}  Assume that $y_1+y_2+y_3\le 6.106$
and that $\rad\ge1.41$.  Then $\sol\ge 0.87$.

{\bf Calculation 9.17.3.4.}  Assume that $y_1+y_2+y_3\le 6.064$
and that $\rad\ge1.41$.  Then $\sol\ge 0.9$.

{\bf Calculation 9.17.3.5.}  Assume that $y_1+y_2+y_3\le 6.032$
and that $\rad\ge1.41$.  Then $\sol\ge 0.91882$.

In the interval arithmetic verification of Calculations
9.17.3, we may assume that $\rad=1.41$ and
that $y_1+y_2+y_3$ is equal to the given upper bound
6.3, 6.192, etc.  To see this, we note that the circumradius
constraint is preserved by a deformation of $S$ that
increases $y_1$, $y_2$, or $y_3$ while keeping fixed
the spherical triangle on the unit sphere at the origin
cut out by $S$.
 We increase $y_1$, $y_2$,
and $y_3$ in this way
until the sum equals the given upper bound.
Then fixing $y_1$, $y_2$, $y_3$, and
one of $y_4$, $y_5$, and $y_6$, we decrease
the other two edges
in such a way as to decrease the solid angle
and circumradius until $\rad=1.41$.

This deformation argument would break down if we encountered
a configuration in which two of $y_4$, $y_5$, and $y_6$ equal
2, but this cannot happen when $\rad(S)\ge 1.41$ because this
constraint on the edges would lead to the contradiction
$$1.41\le\rad(S)\le \rad(S(2.3,2.3,2.3,2,2,2.51))<1.39.$$

It is possible to reduce Calculations 9.17.3 further to the
four-dimensional
situation where two of $y_4$, $y_5$ and $y_6$ are either $2$
or $2.51$.  Consider a simplex $S$ with vertices $0$, $v_1$,
$v_2$, and $v_3$.  Let $p_i$ be the corresponding vertices of
the spherical triangle cut out by $S$ on the unit sphere at the
origin.  Fix the origin, $v_2$, and $v_3$, and vary the vertex $v_1$.
The locus on the unit sphere described as $p_1$
traces out spherical triangles of fixed area is an arc of a
Lexell circle $\Cal C$.  Define the ``interior'' of ${\Cal C}$
to be the points on the side
of $\Cal C$ corresponding to spherical
triangles of smaller area.  The locus 
traced by $v_1$ on the circumsphere (of $S$)
with $|v_1|$ constant
is a circle.
Let ${\Cal C'}$, also a circle, be the radial projection 
of this locus to the unit
sphere.  Define the ``interior'' of $\Cal C'$ be the points coming
from larger $|v_1|$.  The two circles $\Cal C$ and $\Cal C'$
meet at $p_1$, either tangentially or transversely.  The interior
of $\Cal C$ cannot be contained in the interior of $\Cal C'$
because the Lexell arc contains $p_2^*$, the point on the unit
sphere antipodal to $p_2$
\ \cite{FT,p.23}, but the interior of of $\Cal C'$ does
not.  
Furthermore, if $|v_1|=2$ and $|v_2|+|v_3|>4$, then $v_2$ or
$v_3$ lies in the interior of
${\Cal C}$ and $\Cal C'$, so that the circles
have interior points in common.
This means that one can always move $v_1$
in such a way that the solid angle is decreasing,
the circumradius is constant, and the length $|v_1|$ is decreasing
(or constant if $|v_1|=2$).
If any two of $y_4$, $y_5$, and $y_6$ are not at an extreme point,
this argument can be applied to $v_1$, $v_2$, or $v_3$ to
decrease the solid angle.  This proves the reduction.

\smallskip
We are now in a position to prove the lemma for simplices
satisfying $y_1+y_2+y_3\le 6.3$.  Calculation 9.17.3.1 allows
us to assume $\sol\ge 0.767$.  
As in Section 8.6, let $S_y = S(2,2,2,y,y,y)$.
We will rely on the fact that  $\vor(S_y)$ is decreasing in 
$y$, for $y\in [2.26,2.41]$.
In fact, the results of Section 8.6 specialize to the formula
$$\vor(S_y) = {-8\doct y^2\over (12-y^2)^{1/2}(16-y^2)} +
{8\over 3} \arctan\left( {(12-y^2)^{1/2}y^2\over 64-6 y^2}\right),
\tag 9.17.4$$
and the sign of its derivative is determined
by a routine {\it Mathematica\/} calculation.
It is clear that $\sol(S_y)$
is continuous and increasing in $y$.  Since $\sol(S_{2.26})< 0.767$ and
$\sol(S_{2.41})> 0.91882$, our conditions imply that $\sol(S)=
\sol(S_y)$ for some $y\in [2.26,2.41]$.
Let $r(a)$  and $\Vol(R(a,b,c))$ be the functions introduced in
Section 8.6.4.

This suggests the following procedure. Pick
$y$ so that  $\sol(S)\ge \sol(S_y)$.
Calculate
the smallest 
(or at least a reasonably small) $a$ for which 
$$\xi(y,a):=\vor(S_y)- 4\doct
(1-1/a^3) 2r(a)$$
is less than $-1.8\,\pt$.  
Monotonicity (Calculation 9.20.1) and Lemma 8.6.5 imply
that 
$\vor(S)<-1.8\,\pt$, if $y_1+y_2+y_3\ge2(2+a)$.  
To treat the case that remains $(y_1+y_2+y_3<2(2+a)$),
use Calculations 9.17.3 to obtain a new lower bound
on $\sol(S)$, and hence a new value for $y$.
 The procedure is repeated until
$a = 0.016$.  
Calculation 9.17.3.5 completes the argument by
covering the case $y_1+y_2+y_3\le 6.032$.
We leave it to
the reader to check that
$$\xi(2.2626,1.096),\ \xi(2.326,1.053),\ \xi(2.364,1.032),\ 
  \xi(2.391,1.016)$$
are less than $-1.8\,\pt$ and that
$$\sol(S_{2.2626})\le 0.767,\ \sol(S_{2.326})\le 0.83,\ 
\sol(S_{2.364})\le 0.87,\ \sol(S_{2.391})\le 0.9.$$
This completes the proof of Lemma 9.17.\qed

\bigskip
{\bf Lemma 9.18.}  
If
$\rad(S)\ge 1.41$, then
$\vor < -0.419351\sol + 0.2856354$. 

{\bf Proof:}  We adopt the notation and techniques
of Lemma 9.17.
If $\sol\le 0.918819$, then the result follows from Lemma 9.17.
(The proof of 9.17 is independent of the argument that follows
under that restriction on solid angles.)
Let $f(S) = -0.419351\sol(S) + 0.2856354-\vor(S)$.  We
show that $f(S)$ is positive.  We will use the inequality
$f(S)\ge f(S_{\tan})+4\doct\Vol(\hat S\backslash \hat S_{\tan})$
of Lemma 8.6.4.  A routine calculation based on the
formula 9.17.4 shows that
$f(S_y)$ is increasing, for $y\in [2.4085,2.51]$.

If $\sol\ge 0.951385$, then
we appeal to Fejes T\'oth's
convexity argument described in 8.6.  To justify the use
of his argument, we must verify that the cone over $S_{\tan}$
contains the circumcenter of $S$.  The first, second, and
third edges have length 2, and the fourth, fifth, and
sixth edges are between $2.21$ and $2.51$ (Calculation 9.18.2).
These are stronger restrictions on the edges than in
the proof of Lemma 8.6.5, so the justification there applies here
as well.
We 
observe that $\sol(S_{2.4366})< 0.951385$ 
so that
$$f(S)\ge f(S_{\tan}) \ge f(S_y) \ge f(S_{2.4366})\approx 0.000024>0,$$
where $y$ satisfies $\sol(S_y)=\sol(S_{\tan})$.
We may now assume that $0.918819\le \sol\le 0.951385$.

{\bf Lemma 9.18.1.}  The combined volume of the two Rogers
simplices along a common edge of a quasi-regular tetrahedron
$S$ is at least $0.01$.

{\bf Proof:}  The combined volume is at least that
of the right-circular cone of height $a$ and base a 
wedge of
radius $\sqrt{b^2-a^2}$ and dihedral angle  $\dihmin$,
where $a$ is the half-length of an edge and $b$ is a lower
bound on the circumradius of a face with an edge $2a$.
This gives the lower bound of $$
{\pi\over 3}(b^2-a^2)a {\dihmin\over 2\pi}>0.1439
(b^2-a^2).$$  
We minimize $b$ by setting $b^2=\eta(2,2,2a)^2=4/(4-a^2)$.
Then $b^2-a^2= (a^2-2)^2/(4-a^2)$, which is decreasing in $a\in[1,2.51/2]$.
  Thus, we obtain a lower bound on $b^2-a^2$ by setting
$a=2.51/2$, and this gives the estimate of the lemma.\qed

We remark that $\sol(S_{2.4085})<0.918819$.
If $1.15\le (y_1+y_2+y_3-4)/2$, then Lemma 9.18.1
and Section 8.6 give
$$f(S)\ge f(S_{2.4085}) + 4\doct (1-1/1.15^3) 0.01> 0.$$
(Analytic continuation is not required here, because of the constraints
on the edges in Calculation 9.18.2.)
We may now assume that $y_1+y_2+y_3\le 6.3$.
To continue, we
need a few more calculations.

{\bf Calculation 9.18.2.}  If $\sol\ge 0.918$, then
$y_4, y_5, y_6\ge 2.21$. 

{\bf Calculation 9.18.3.1.} If $y_1+y_2+y_3\le 6.02$ and $\rad\ge1.41$,
then $\sol\ge 0.928$. 

{\bf Calculation 9.18.3.2.} If $y_1+y_2+y_3\le 6.0084$ and $\rad\ge1.41$,
then $\sol\ge 0.933$.  

{\bf Calculation 9.18.3.3.}  If $y_1+y_2+y_3\le 6.00644$ and $\rad\ge1.41$,
then $\sol\ge 0.942$. 

In the verification of these calculations, we make the same
reductions as in the Calculations 9.17.3.

Adapting the procedure of
Lemma 9.17 and Lemma 8.6.5, 
we find that a lower bound on the solid angle leads
to an estimate of a constant $a$ with the property that 
$f(S)>0$ whenever $y_1+y_2+y_3\ge 4+2a$.  
That is, we pick $a$ so that
$\xi_1(y,a):=f(S_y)+8\doct(1-1/a^3)r(a)$ is positive, where
$y$ is chosen so that $\sol(S_y)$ is a lower bound
on the solid angle.  
The values
$$\xi_1(2.4085,1.01),\ \xi_1(2.4165,1.0042),\ \xi_1(2.42086,1.00322),\ 
\xi_1(2.4286,1.0017)$$
are all positive.
This yields the bound
$y_1+y_2+y_3\le 6.0034$.

Assume that $S$ satisfies $y_1+y_2+y_3\le 6.0034$ and
$\rad(S)\ge 1.41$.  Then by Calculation 9.18.3.3,
$\sol(S)=\sol(S_{\tan})\ge 0.942$. Also
$\rad(S_{\tan}) > \rad(S)/1.0017\ge1.41/1.0017 > 1.4076$, because
rescaling $S_{\tan}$ by a factor of $1.0017$ gives
a simplex containing  $S$.  This means that $S_{\tan}$
satisfies the hypotheses of the following calculation.
Calculation 9.18.4 completes the proof of
Lemma 9.18. \qed

{\bf Calculation 9.18.4.}  If $\rad(S_{\tan})\ge 1.4076$
and $0.933\le \sol(S_{\tan})\le 0.951385$, then $f(S_{\tan})>0$.

In this verification, we may assume that the fourth
fifth, and sixth edges of $S_{\tan}$ are at least $2.27$.
for otherwise the circumradius is at most 
$$\rad(S(2,2,2,2.27,2.51,2.51))
< 1.4076.$$
In the verification of 9.18.4, we also rely on the fact that
$$f_1(x_4,x_5,x_6):= f(S(2,2,2,\sqrt{x_4},\sqrt{x_5},\sqrt{x_6}))$$
is increasing in $(x_4,x_5,x_6)\in [2.27^2,2.51^2]$.  Here is
a sketch justifying this fact.  The details were carried out
in {\it Mathematica\/} with high-precision arithmetic.
The explicit formulas of Section 8.6 lead to an expression for
$\partial f_1/\partial x_4$ as
$${W(x_4,x_5,x_6)(-x_4+x_5+x_6)\over (-16+x_4)^2 \Delta^{3/2}},$$
where $W$ is a polynomial in $x_4$, $x_5$, and $x_6$
(with 13 terms).  To show that
$W$ is positive, expand it in a 
Taylor polynomial about $(x_4,x_5,x_6)=(2.27^2,2.27^2,2.27^2)$
and check that the lower bound of Inequality 7.1 is positive.

\bigskip

{\bf Calculation 9.19.}  If $y_4\in[2,2.1]$, then $\sol<0.906$.

For $a\le b\le c$, we have $\Vol(R(a,b,c)) 
= a((b^2-a^2)(c^2-b^2))^{1/2}/6$.
 The final four calculations
are particularly simple (to the extent that
any of these calculations are simple), 
since they involve a single variable.
They were verified in {\it Mathematica\/} with rational arithmetic.

{\bf Calculation 9.20.1.}  The function $(1-1/a^3)\Vol(R(a,\eta(2,2,2a),1.41))$
is increasing on $[1,1.15]$. 

{\bf Calculation 9.20.2.}  The function $\Vol(R(a,\eta(2,2,2a),1.41))$
is decreasing for $a\in [1,1.15]$.

{\bf Calculation 9.20.3.}  For $a\in[1,1.15]$ we have
$$\Vol(R(a,\eta(2,2,2a),1.41))< 
\Vol(R(a,\eta(2,2.51,2a),1.41)).$$

{\bf Calculation 9.21.}  Let $p_1=p_1(y_1,\ldots,y_6)$ 
be the point in Euclidean
space 
introduced in Section 8.6.7. For $y\in [2.3,2.51]$, 
$${d\phantom{d}\over dy} (\eta(2.51,2.51,y)-y^2/4)^{1/2}
  > -0.75 >  {d\phantom{d}\over dy} (|p_1(y,2.51,2,2,2,2.51)|^2
   - y^2/4)^{1/2}.$$

\vfill\eject
\bigskip
\Refs

[FT] L. Fejes T\'oth, Lagerungen in der Ebene, auf der Kugel
	und im Raum, Springer-Verlag, 1953.

[Int] G\"otz Alefeld and J\"urgen Herzeberger, Introduction to
        Interval Computations, Academic Press, 1983.

[H1] Thomas C. Hales, Remarks on the Density of Sphere Packings
        in Three Dimensions, Combinatorica, 13 (2) (1993), 181--197.

[H2] Thomas C. Hales, The sphere packing problem, Journal of Computational
and Applied Math, 44 (1992) 41--76.

[H3] Thomas C. Hales, The status of the Kepler conjecture, 
	Math. Intelligencer, (1994).

[H4] Thomas C. Hales, Sphere Packings II, to appear in 
	{\it Discrete and Computational Geometry}.

[H5] Thomas C. Hales, Sphere Packings III$_\alpha$, preprint.

[H6] Thomas C. Hales, Packings,\hfill\break\quad
{\tt http://www-personal.math.lsa.umich.edu/\~\relax hales/packings.html}

[CA] John F. Hart et. al, Computer Approximations, John Wiley
and Sons, 1968.

[IEEE]  IEEE Standard for Binary Floating-Point Arithmetic, ANSI/IEEE
                Std. 754--1985.

[NR] William H. Press et al., Numerical Recipes in C,
        Chapter 20. Less-Numerical Algorithms, Second Edition,
        Cambridge University Press, 1992.

[R] C. A. Rogers, The Packing of Equal Spheres, Proc. London Math. Soc.
        (3) 8 (1958), 609--620.

[W] What Every Computer Scientist Should Know about Floating-Point
        Arithmetic, Computer Surveys, March 1991, Association for
        Computing Machinery Inc.

\endRefs

\vfill\eject

\centerline{\bf Appendix.  Proof of Theorem 6.1}
\medskip
\centerline{Douglas J. Muder}
\bigskip


\subheading{Notation and observations} Let $P$ be a point
of degree d.  If we consider $P$ to be the center of the
configuration, then the {\it first rim} of points will be the
set of $d$ points adjacent to $P$, and the {\it second rim}
will be those at distance 2 from $P$. Let $\delta(P)$ be
the sum of the degrees of the first rim points. If $d=6$, it
is easy to see that the number of points on the second rim
is $s = \delta(P) - 24$, and the total number of points of
distance at most 2 from $P$ is $\delta(P) - 17$. The
second rim is thus an $s$-gon. This $s$-gon and the
triangulation of the $P$-side of the $s$-gon will be
referred to as the {\it inner graph}. The inner-graph
degree of a second rim point will be called its {\it inner
degree}. The number of second-rim points with inner
degree 5 is equal to the number of degree 4 points in the
first rim, and the number of second-rim points with inner
degree 3 is equal to the number of degree 6 points on the
first rim. All other second-rim points have inner degree 4.
Points beyond the second rim will be called {\it extra
points}. If there are no extra points, then there are at least
2 second-rim points whose degrees are equal to their inner
degrees. Let $N^i_{\Delta}(P)$ be the number of points of
degree $\Delta$ in rim $i$ around $P$. Notice that if
$d=6$, then
$$\delta(P) = 30 - N^1_4(P) + N^1_6(P).$$
Further, Euler's formula gives $N-12 = N_6 - N_4$. Putting
this together with 5.1.1 and 5.1.5, we see
$$\gather
13 \leq 12 + N_6 - N_4\leq 15 \\
1 \leq 1 + N_4 \leq N_6 \leq 3 + N_4 \leq 5.
\endgather
$$

\proclaim{Lemma 1} There is no triangle whose
vertices all have degree six. \endproclaim

\proclaim{Lemma 2} Suppose $N_6\geq 3$. Let $\bold
G_6$ be the subgraph of points of degree six.  Either
$\bold G_6$ is three points in two components, or it is a
single (open or closed) path, with no other edges.
\endproclaim
 
With these lemmas and this notation, we consider the
different possible values of $N_6$.
 
\subheading{$N_6 = 1$} This forces $N_4=0$. Let $P$ be the
degree 6 point. Then $\delta(P) = 30$, all 13 points are in the
inner graph, and all 6 second-rim points have inner degree
4. But at least two of these second-rim points have inner
degree equal to degree. Contradiction.
 
\subheading{$N_6 = 2$} If the two degree 6 points are
adjacent, let $P$ be either one. Now $\delta(P) = 31 -
N_4^1$, and we see that there are
no extra points, and all second rim points have degree 5.
All but $N_4^1$ second rim points have inner degree less
than 5, and at least two of them will have inner degree
equal to degree.
So $N_4\ge N_4^1 = 2$, and $N\le 12$. Contradiction.

If the degree six points are not adjacent, then the $(6,6)$
forces $N\geq 14$, so $N_4 = 0$.  If $\delta(P)$ is either
of the two degree six points, $\delta(P) = 30$, so the
second rim has 6 points, all of inner degree 4, and there is
one extra point. If the extra point has degree 5, it is
adjacent to all but one of the second rim points. This
unique second-rim point is then part of a quadrilateral,
which can only be triangulated by the diagonal edge that
does not include the extra point. This creates two
second-rim points of degree 6, in violation of our
assumptions. Therefore the extra point has degree 6, and
we have 6.2.
 
\subheading{$N_6 = 3$} Either $\bold G_6$ is an open
path or it has two components. In either case there is a
$(6,6)$, so $N\geq 14$ and $N_4=0$ or $1$.
 
In the first case let $P$ be the center of the path. Now
$\delta(P) = 32 -N^1_4$. We see that there are no extra
points and any degree 4 point is in the first rim. The
second rim has at most one point of inner degree greater than
4, and two points whose inner degrees are equal to their
degrees. So there must be a point of degree 4 on the second
rim. Contradiction.
 
In the second case 5.1.8 forces $N=15$, hence $N_4=0$.
If $P$ is one of the points in the two-point component of
$\bold G_6$, then $\delta(P)=31$, so there is one extra
point and 7 second-rim points. Either one or two of the
second rim points is not connected to the extra point. In
either case at least one of these two has degree = inner
degree, which is at most 4. Contradiction.
 
\subheading{$N_6 = 4$} Now $\bold G_6$ is either an
open or closed path of length 4. If open, then by 5.1.8,
$N=15$ and $N_4=1$. If $P$ is either of the interior
degree 6 points, then $\delta(P)=32 - N_4^1$, the second
rim is an $(8-N^1_4)$-gon with $N_4^1$ points of
inner degree at least 5 and $N_4^1$ extra points. We
cannot have $N_4^1=0$, or else there would be at least two
second rim points of degree 4, so $N_4^1 = 1$. This
holds no matter which interior degree 6 point we started
with, so the degree 4 point must be adjacent to both. It also
must be adjacent to one of the other two, or else there is a
$(6,6,4)$.  Therefore there is a degree 6 point whose
first rim degree sequence is $6/4/6/5/5/5$, producing a
second rim inner degree sequence of $4/3/5/3/4/4/4$.
The point of inner degree 5 is the fourth neighbor of the
degree four point, and cannot be degree 6 without making
$\bold G_6$ closed. Therefore there is an edge connecting
the two points of inner degree 3. Adding this edge
to the inner graph creates a hexagon with all points of
``inner'' degree 4, and the extra point in its interior. The
extra point is connected to 5 of these points, and the sixth
has degree 4. Contradiction.
 
If $\bold G_6$ is a closed path, then a degree 4
point $Q$ must triangulate this quadrilateral. Using $Q$
as the center, there are 4 first rim points and 8 second
rim points. The second rim can be drawn in a square with the
4 points of inner degree 3 on the corners and the 4 points of
inner degree 4 at the midpoints. There is at most one
additional degree 4 point (other than $Q$). If it is
anywhere but at one of the second-rim midpoints, a
$(6,6,4)$ exists. The four midpoints cannot be adjacent to
any additional second-rim points without creating an
illegal triangle or quadrilateral. Since there can be at most
one more degree 4 point, at least 3 of the 4 second-rim
midpoints must be joined
to the extra point(s). Two midpoints of consecutive sides
can't be joined to the same extra point without forcing the
corresponding corner to be degree 4, which is impossible.
So there must be two extra points, and thus no degree four
points other than $Q$. One extra point must be joined to
the midpoints of each pair of parallel sides of the square.
But all these edges cannot be drawn without intersecting.
 
\subheading{$N_6 = 5$} Now $\bold G_6$ is a path of
length 5, closed so that a $(6,6,6)$ doesn't exist. Any
point in this pentagon is nonadjacent to two others, so the
situation of 5.1.8 applies in 5 different ways. Let $P, Q, R\in
\bold G_6$ with $P$ and $Q$ adjacent, but neither adjacent
to $R$. Then there are 8 points adjacent to either $P$ or $Q$,
and 6 points adjacent to $R$, so there are exactly two
common points in these two sets. In our case these two
points are precisely the other two points of degree 6. This
means that no point inside or outside the pentagon can be
connected to more than two points of the pentagon. Thus
the only possible configuration for the 15 points is to form
three concentric pentagons, with $\bold G_6$ as the
middle one. Prior to triangulating the inner and outer
pentagons, all of their points have degree 4. Triangulating
either creates more degree 6 points. Contradiction.
 
\demo{Proof of Lemma 1} Let $P_1$, $P_2$, and
$P_3$ be the vertices of such a triangle. If $P_j$ has no
neighbors of degree 4, then the only possibility is
$\delta(P_j)=32$, which forces $N^1_5=4, N^1_6=2$, and
$N=15$. All second-rim points have inner degree at most 4,
and there are no extra points. Therefore two second-rim
points have degree and inner degree 4. Call them $Q_1$ and
$Q_2$. Therefore, $N_4\geq 2$, $N_6\geq 5$, and
$N_6^2\geq 2$. Let $R_1$ and $R_2$ be second-rim points
of degree 6. Each $Q_k$ and $R_{\ell}$ must be adjacent, or
else $(P_j, R_{\ell}, Q_k)$ is a $(6,6,4)$. But all edges from
$Q_j$ are inner-graph edges. So $Q_1R_1Q_2R_2$ must be
a second-rim quadrilateral, which is impossible. Thus
$N_4^1(P_j)\geq 1$ for each $P_j$.
 
This forces $N_4=2$. At least one of the degree 4 points is
adjacent to two of the $P_j$. One of these two, say $P_1$, is
not adjacent to the other degree 4 point. Since $\delta(P)\leq
32$, we are left with four possible first-rim degree sequences
for $P_1$: (1) 6/6/4/5/5/5; (2) 6/6/4/6/5/5;
(3) 6/6/4/5/6/5; or (4) 6/6/4/5/5/6.
 
The last three possibilities are easily dealt with.
Sequence (3) contains a $(6,6,4)$.  Sequence (4) creates
two 6-6-6 triangles joined at an edge. Let $P_1$ and
$P_2$ be the vertices of the edge, and $P_3$ and $P_4$
be the other vertices of the triangles.  The conditions
$N_4^1(P_j)\geq 1$ and $N_4=2$ force $P_3$ and $P_4$
to have sequence (1). In (2) we can invoke 5.1.8 with the 3
degree 6 points in the first rim. Call them $R_1$, $R_2$ and
$T$, with $R_1$ adjacent to $R_2$. In
order for the numbers to work out, there can be precisely
two points which are adjacent to $T$ as well as one of
the $R_j$. But in (2) the center,
the first rim degree 4 point, and the second rim neighbor of
the first rim degree 4 point must all fit this description.
 
Sequence (1) is more difficult to eliminate. Let $R_1,
\dots R_6$ be the first rim points listed in the order of
(1). Let $S_1, \dots, S_7$ be the second rim points, with
$S_1$ being adjacent to both $R_1$ and $R_2$, and
$S_7$ having only $R_1$ as a first rim neighbor. There is
at least one degree 6 point in the second rim---if it is at
$S_1$, then $\delta(R_1)=32$; if it is anywhere else we
can invoke 5.1.8.  In either case, $N=15$ and there are
two degree 6 points outside the first rim.  These must lie in
$\{S_6, S_7, S_1, S_3\}$ to avoid making a $(6,6,4)$
with $R_1$ and $R_3$.  At most
one of these points can be adjacent to
$R_1$ (since $\delta(R_1)<33$), so $S_3$ must
have degree 6. Now $S_3$ must connect  to both the
degree 6 and degree 4 points in $\{S_6, S_7, S_1\}$ to
avoid making a $(6,6,6)$ or $(6,6,4)$ with the center
point. But $S_3$ has inner degree 5 and so can connect to
at most one of those points. Contradiction.
\qed\enddemo
 
\demo{Proof of Lemma 2}
If $\bold G_6$ has at least three points and no triangles,
then there exists a $(6,6)$ forcing $N\geq 14$. Any point
in $\bold G_6$ can be adjacent to at most two
others, or else either a triangle or a $(6,6,6)$ is created.
Therefore each component is a path.
There cannot be three or more components without
producing a $(6,6,6)$. If $\bold G_6$ has two components,
each of them must be a complete graph for the same
reason. Since there are no triangles, both of the
components must have at most two points. If there are two
components with two points each, then 5.1.8 forces $N=15$.
Each component is adjacent to 8 non-$\bold G_6$ points of
the original graph, and not adjacent to 3. Now 5.1.8 forces
each point in the other component to be adjacent to these 3.
But two points can have only two common neighbors.
So the only two-component $\bold G_6$ is the
one described in the lemma.
\qed
\enddemo

\enddocument\bye